\documentclass{article}
\usepackage{amssymb}
\usepackage{amsfonts}
\usepackage{amsmath}

\newtheorem{theorem}{Theorem}

\newtheorem{corollary}[theorem]{Corollary}

\newtheorem{definition}[theorem]{Definition}

\newtheorem{lemma}[theorem]{Lemma}

\newtheorem{proposition}[theorem]{Proposition}

\input{tcilatex}
\begin{document}

\title{Characteristic classes as obstructions\\
to local homogeneity}
\author{Erc\"{u}ment Orta\c{c}gil}
\maketitle

This note is an extended version of my presentation at the "Focused Research
Workshop on Exterior Differential Systems and Lie Theory" held at the Fields
Institute in Toronto, December 9-13, 2013. Some ideas and constructions of
this note crystallized during my communication with Anthony D. Blaom. I am
deeply grateful to him and also to P.J. Olver for his constant help and
encouragement.

\bigskip

\textbf{Contents}

\bigskip

1) Summary

2) Prehomogeneous geometries

3) The algebroid of a prehomogeneous geometry

4) Curvature

5) Cartan algebroids

6) Characteristic classes on the base

7) Higher order characteristic classes

8) Dependence on the isomorphism class

9) A. Cartan connections

10) B. Chern-Simons forms

11) C. Uniformization number and representations

12) D. The adjoint representation

\section{Summary}

Let $M$ be a smooth manifold with $\dim M=n$ and $\left\{ U_{\alpha
}\right\} $ be an atlas with transition functions $\phi _{\beta }\circ \phi
_{\alpha }^{-1}:U_{\alpha }\cap U_{\beta }\rightarrow \mathbb{R}^{n}.$ Are
there any "global invariants" of $M$ (at least for certain $M)$ which depend
on the $k$-th order derivatives of $\phi _{\beta }\circ \phi _{\alpha }^{-1}$
for arbitrarily large $k$ (as $M$ varies)? Equivalently, do the higher order
derivatives play any role in global differential geometry? This note is the
outcome of our efforts over a period of 20 years and gives, we hope, an
affirmative answer to this question in Section 7. Our method produces also
obstructions to $k$-flatness as defined in [9]. We will shortly outline here
the construction of these invariants which turn out to be the "old friends"
but seen with a new eye.

A prehomogeneous geometry $\varepsilon \mathcal{G}_{k}$ $(phg$ for short) of
order $k$ on $M$ is a very special transitive Lie groupoid on $M.$ The
integer $k\geq 0$ indicates the order of jets involved in the definition of $%
\varepsilon \mathcal{G}_{k}.$ Now $\varepsilon \mathcal{G}_{0}$ is an
absolute parallelism on $M,$ $k=1$ for Riemannian geometry but can be
arbitrarily large for parabolic geometries (like projective and conformal
geometries) \textit{as defined }in Section 2. The curvature $\mathcal{R}_{k}$
of $\varepsilon \mathcal{G}_{k}$ vanishes if and only if the $PDE$ defined
by $\varepsilon \mathcal{G}_{k}$ is locally solvable. In geometric terms,
this is equivalent to the local homogeneity of $M$ in the way imposed by $%
\varepsilon \mathcal{G}_{k}.$ With the assumption of completeness and simple
connectedness, $M$ becomes a globally homogeneous space $G/H$ possibly with
noncompact $H.$ In fact, compactness of $H$ forces $k\leq 1$.\ We have $%
k\leq $ $\dim Nil(\mathfrak{h})+1$ so that $k\leq 1$ also if $H$ is
semisimple. However $\varepsilon \mathcal{G}_{k}$ is not modeled on some
fixed $G/H$ chosen beforehand.

The algebroid $\mathfrak{\varepsilon G}_{k}\rightarrow M$ of $\varepsilon 
\mathcal{G}_{k}$ is a very special vector bundle filtered by jets. The
Chern-Weil construction applied to the curvature $\mathfrak{R}_{k}$ of $%
\varepsilon \mathfrak{G}_{k}\rightarrow M$ establishes the Pontryagin
algebra $\mathcal{P}^{\ast }(M,\varepsilon \mathfrak{G}_{k})\subset
H_{dR}^{\ast }(M,\mathbb{R)}$ as an obstruction to local homogeneity. In
other words, the well known characteristic classes of vector bundles become
obstructions to integrability once they are restricted to this particular
subset of vector bundles. These obstruction depend on first order jets and
are topological. Using the projections $\mathcal{G}_{k}\rightarrow \mathcal{G%
}_{r},$ $0\leq r\leq k,$ we define the higher order Pontryagin algebras $%
\mathcal{P}^{\ast }(\mathcal{G}_{r}^{\bullet },\varepsilon \mathfrak{G}%
_{k})\subset H^{\ast }(\mathcal{G}_{r}^{\bullet },\mathbb{R})$ where $%
\mathcal{G}_{r}^{\bullet }\rightarrow M$ is the principal bundle of the
groupoid $\mathcal{G}_{r}$ so that $\mathcal{G}_{0}=M\times M,$ $\mathcal{G}%
_{0}^{\bullet }=M.$ For $1\leq r\leq k$, $\mathcal{P}^{\ast }(\mathcal{G}%
_{r}^{\bullet },\varepsilon \mathfrak{G}_{k})$ is trivial as a subalgebra of 
$H_{dR}^{\ast }(\mathcal{G}_{r}^{\bullet },\mathbb{R})$ due to the
contractibility of the fibers of $\mathcal{G}_{k}^{\bullet }\rightarrow 
\mathcal{G}_{r}^{\bullet }.$ However, the representatives of $\mathcal{P}%
^{\ast }(\mathcal{G}_{r}^{\bullet },\varepsilon \mathfrak{G}_{k})$ are right
invariant forms on the principal bundle $\mathcal{G}_{r}^{\bullet
}\rightarrow M$ and generate a subalgebra $\widehat{\mathcal{P}^{\ast }}(%
\mathcal{G}_{r}^{\bullet },M)$ in the subcomplex of right invariant forms
which computes the algebroid cohomology of $\mathcal{G}_{r}^{\bullet
}\rightarrow M.$ The subalgebra $\widehat{\mathcal{P}^{\ast }}(\mathcal{G}%
_{r}^{\bullet },M)\subset H_{inv}^{\ast }(\mathcal{G}_{r}^{\bullet
},M)=H^{\ast }(M,\mathfrak{G}_{r}),$ we believe, need not be trivial for $%
1\leq r\leq k$ and gives obstructions to local homogeneity which depend on
jets of order $k$. All these obstructions depend on the isomorphism class $%
\left[ \varepsilon \mathcal{G}_{k}\right] $ of $\varepsilon \mathcal{G}_{k}.$
In view of the definition of $\left[ \varepsilon \mathcal{G}_{k}\right] ,$
the assignment $\left[ \varepsilon \mathcal{G}_{k}\right] \Rightarrow 
\mathcal{P}^{\ast }(\mathcal{G}_{r}^{\bullet },\varepsilon \mathfrak{G}_{k})$
is tantamount to the assignment of certain invariants to the moduli space of
connections on the principal bundle $\mathcal{G}_{k}^{\bullet }\rightarrow M$
as in gauge theory. In Section 8 we interpret the above Pontryagin algebras
as obstructions to the existence of certain Cartan connections. In Section 9
we show that the Chern-Simons forms arise naturally in the present framework
but with a surprisingly different interpretation.

\section{Prehomogeneous geometries}

Let $M$ be a smooth manifold with $\dim M=n\geq 2$ and $j_{k}(f)^{p,q}$ be
the $k$-jet of the local diffeomorphism $f$ with source at $p$ and target at 
$q.$ We call $j_{k}(f)^{p,q}$ a $k$-arrow from $p$ to $q.$ Clearly $%
j_{0}(f)^{p,q}=(p,q).$ Let $\mathcal{U}_{k}^{p,q}$ denote the set of all $k$%
-arrows from $p$ to $q.$ With the composition and inversion of $k$-arrows,
the set $\mathcal{U}_{k}\overset{def}{=}\cup _{p,q\in M}\mathcal{U}%
_{k}^{p,q} $ of all $k$-arrows on $M$ has the structure of a groupoid which
we call the universal $k$-th order groupoid on $M$. The set $\mathcal{U}%
_{k}^{p,p},$ $p\in M,$ is a Lie group and called the vertex group of $%
\mathcal{U}_{k}$ at $p$. A choice of coordinates around $p$ identifies $%
\mathcal{U}_{k}^{p,p}$ with the $k$-th order jet group $G_{k}(n)$ in $n$
variables. We define $G_{0}(n)$ as the set containing one point. The
projection of jets induces a homomorphism $\pi _{k+1,k}:\mathcal{U}%
_{k+1}\rightarrow \mathcal{U}_{k}$ of groupoids and we have the sequence of
projections

\begin{equation}
...\rightarrow \mathcal{U}_{k+1}\rightarrow \mathcal{U}_{k}\rightarrow
...\rightarrow \mathcal{U}_{1}\rightarrow \mathcal{U}_{0}=M\times M
\end{equation}

Note that (1) can be restricted to the vertex groups at $p.$ The set $%
\mathcal{U}_{k}^{e,\bullet }\overset{def}{=}\cup _{x\in M}\mathcal{U}%
_{k}^{e,x}$ is a principal bundle with the structure group $\mathcal{U}%
_{k}^{e,e}\simeq G_{k}(n)$ where $e\in M$ is some basepoint and (1)
restricts also to these principal bundles.

In this note we will be interested in certain subgroupoids $\mathcal{G}%
_{k}\subset \mathcal{U}_{k}.$ For $s\leq k,$ we denote $\pi _{k,s}\mathcal{G}%
_{k}$ by $\mathcal{G}_{s}\subset \mathcal{U}_{s}.$

\begin{definition}
A prehomogeneous geometry ($phg)$ of order $k$ on $M$ is a subgroupoid $%
\mathcal{G}_{k+1}\subset \mathcal{U}_{k+1}$ satisfying

$i)$ $\mathcal{G}_{0}=\mathcal{U}_{0}=M\times M$

$ii)$ $\mathcal{G}_{k+1}\simeq \mathcal{G}_{k}$ and $k$ is the smallest such
integer
\end{definition}

So $\mathcal{G}_{k}\subset \mathcal{U}_{k}$ is an imbedded submanifold
consisting of certain $k$-arrows of $\mathcal{U}_{k}$ closed under
composition and inversion. $i)$ states that $\mathcal{G}_{k}$ is transitive
on $M$, i.e., for any $p,q\in M$ there exists a $k$-arrow of $\mathcal{G}%
_{k} $ from $p$ to $q.$ Let $\varepsilon :\mathcal{G}_{k}\rightarrow 
\mathcal{G}_{k+1}$ denote the inverse of the isomorphism given by $ii)$ so
that $\mathcal{G}_{k+1}=\varepsilon \mathcal{G}_{k}$. Now $ii)$ states that
above any $k$-arrow $j_{k}(f)^{p,q}$ in $\mathcal{G}_{k}^{p,q}$ there exists
a unique $\left( k+1\right) $-arrow (namely $\varepsilon j_{k}(f)^{p,q})$
and this 1-1 correspondence preserves composition and inversion of arrows
since $\varepsilon $ is an isomorphism of groupoids. The second condition in 
$ii)$ states that $\mathcal{G}_{r}$ projects onto $\mathcal{G}_{s}$ with
nontrivial kernel for $1\leq s+1\leq r\leq k$. Since $\mathcal{G}_{k}$ is
transitive, this condition holds if and only if it holds at one (hence all)
vertex group. Note that $\varepsilon $ restricts to the vertex groups $%
\mathcal{G}_{k}^{p,p}$ and also to the principal bundle $\mathcal{G}%
_{k}^{e,\bullet }\rightarrow M.$ As we will see shortly, many (if not all)
geometric structures are particular $phg$'s. We will add a third condition $%
iii)$ to Definition 1 below when it is needed.

Choosing coordinates $(U,x^{i}),$ $(V,y^{i})$, elements of $\mathcal{U}%
_{k+1}^{p,q}$ with $p\in U$, $q\in V$ can be expressed locally as $%
(x^{i},y^{i},f_{j_{1}}^{i},f_{j_{2}j_{1}}^{i},...,f_{j_{k+1}...j_{1}}^{i}).$
Since $\varepsilon \mathcal{G}_{k}\subset \mathcal{U}_{k+1}$ is a
submanifold, locally it is defined by a set of independent equations

\begin{equation}
\Phi ^{\alpha
}(x^{i},y^{i},f_{j_{1}}^{i},f_{j_{2}j_{1}}^{i},...,f_{j_{k+1}...j_{1}}^{i})=0%
\text{ \ \ \ }1\leq \alpha \leq \dim \mathcal{U}_{k+1}-\dim \varepsilon 
\mathcal{G}_{k}
\end{equation}

The functions $\Phi ^{\alpha }$ are surely not unique and the study of their
invariance properties gives rise to a subtle local theory which we will not
touch here. Note that (2) puts no restriction on the variables $x^{i},y^{i}$
by $i)$. Since $f_{j_{k+1}...j_{1}}^{i}$ is determined by $%
x^{i},y^{i},f_{j_{1}}^{i},$...$f_{j_{k}...j_{1}}^{i}$ by $ii),$ we can solve 
$f_{j_{k+1}...j_{1}}^{i}$ in terms of $x^{i},y^{i},f_{j_{1}}^{i},$...$%
f_{j_{k}...j_{1}}^{i}$ and rewrite (2) in an equivalent form.

We now fix an "initial condition" $(\overline{x}^{i},\overline{y}^{i},%
\overline{f}_{j_{1}}^{i},\overline{f}_{j_{2}j_{1}}^{i},...,\overline{f}%
_{j_{k+1}...j_{1}}^{i})$ satisfying (2). We search for a diffeomorphism $%
f:U\rightarrow f(U)\subset V$ such that $(x^{i},f^{i}(x),\frac{\partial
f^{i}(x)}{\partial x^{j_{1}}},$ $...,$ $\frac{\partial ^{k}f^{i}(x)}{%
\partial x^{j_{k+1}}...\partial x^{j_{1}}}$) solves (2) for all $x\in $ $U$
and also satisfies $f^{i}(\overline{x})=\overline{y}^{i},$ $\frac{\partial
f^{i}(\overline{x})}{\partial x^{j_{1}}}=\overline{f}_{j_{1}}^{i},...,\frac{%
\partial ^{k}f^{i}(\overline{x})}{\partial x^{j_{k+1}}...\partial x^{j_{1}}}=%
\overline{f}_{j_{k+1}...j_{1}}^{i}.$ This interpretation shows that $%
\varepsilon \mathcal{G}_{k}$ is a nonlinear $PDE$ of order $k+1$ defined on
the universal pseudogroup $Diff_{loc}(M)$ of all local diffeomorphisms of $M$
and is locally of the form (2). The $(k+1)$-arrows of $\varepsilon \mathcal{G%
}_{k}$ are the initial conditions. In a coordinate free language, let $%
\alpha _{k+1}^{p,q}\in \varepsilon \mathcal{G}_{k}.$ A local diffeomorphism $%
f:U\rightarrow f(U)$ with $f(p)=q$ is a local solution of $\varepsilon 
\mathcal{G}_{k}$ on $U$ satisfying the initial condition $\alpha
_{k+1}^{p,q} $ if

$i)$ $j_{k+1}(f)^{p,q}=\alpha _{k+1}^{p,q}$

$ii)$ $j_{k+1}(f)^{x,f(x)}\in \varepsilon \mathcal{G}_{k}$ for all $x\in U$

A local solution, if it exists, satisfies all its $\left( k+1\right) $%
-arrows as initial conditions.

\begin{proposition}
If $f,g$ are two local solutions with $j_{k+1}(f)^{p,q}=j_{k+1}(g)^{p,q}$,
then $f=g$ on their common domain of definition.
\end{proposition}

Proposition 2 states that local solutions, if they exist, are unique. This
can be seen roughly by noting that $\varepsilon $ expresses jets of order $%
k+1$ in terms of the lower order jets so that the Taylor expansion of a
local solution satisfying some initial condition is determined by this
initial condition.

\begin{definition}
$\varepsilon \mathcal{G}_{k}$ is locally solvable if all its $\left(
k+1\right) $-arrows integrate to local solutions as above.
\end{definition}

Suppose $\varepsilon \mathcal{G}_{k}$ is locally solvable and let $%
\widetilde{\varepsilon \mathcal{G}_{k}}$ denote the set of all local
diffeomorphism obtained by integrating the $\left( k+1\right) $-arrows of $%
\varepsilon \mathcal{G}_{k}.$ Since $\varepsilon \mathcal{G}_{k}$ is a
groupoid, we easily see that $\widetilde{\varepsilon \mathcal{G}_{k}}$ is a
pseudogroup and we say that $\varepsilon \mathcal{G}_{k}$ integrates to $%
\widetilde{\varepsilon \mathcal{G}_{k}}$. Therefore, if $\varepsilon 
\mathcal{G}_{k}$ is locally solvable, then $M$ is locally homogeneous in the
way imposed by $\varepsilon \mathcal{G}_{k}$.

Now suppose that $\varepsilon \mathcal{G}_{k}$ is locally solvable. Let $%
f\in \widetilde{\varepsilon \mathcal{G}_{k}}$ with $f(p)=q$ and $\gamma $ be
a (continuous) path from $\ p$ to some point $r.$ Using Proposition 2 we can
"analytically continue" $j_{k+1}(f)^{p,q}$ along this path but we may not be
able to end up with a $(k+1)$-arrow with source at $r.$

\begin{definition}
If elements of $\widetilde{\varepsilon \mathcal{G}_{k}}$ can be analytically
continued indefinitely along paths, then $\widetilde{\varepsilon \mathcal{G}%
_{k}}$ is complete.
\end{definition}

If $f\in \widetilde{\varepsilon \mathcal{G}_{k}}$ is the restriction of some
(unique!) global transformation $\widetilde{f}\in Diff(M)$, we call $f$
globalizable.

\begin{definition}
$\widetilde{\varepsilon \mathcal{G}_{k}}$ is globalizable if all $f\in 
\widetilde{\varepsilon \mathcal{G}_{k}}$ are globalizable.
\end{definition}

Hence if $\widetilde{\varepsilon \mathcal{G}_{k}}$ is globalizable then we
obtain a global transformation group $G$ which acts transitively and
effectively on $M.$ We call this data a Klein geometry $(G,M)$ which we can
identify with the homogeneous space $G/H=M$ where $H$ is the stabilizer at
some point. Note that this identification is not canonical and depends on
the choice of a base point. Obviously $\widetilde{\varepsilon \mathcal{G}_{k}%
}$ is complete if it is globalizable. Conversely, let $f\in \widetilde{%
\varepsilon \mathcal{G}_{k}}$ with $f(p)=q.$ Assuming that $\widetilde{%
\varepsilon \mathcal{G}_{k}}$ is complete and $M$ is simply connected, we
define a map $\widetilde{f}:M\rightarrow M$ as follows: for any $r\in M$, we
choose a path from $p$ to $r,$ continue $f$ along this path up to $r$ and
define $\widetilde{f}(r)$ to be the value of this continuation. A standard
monodromy argument using simple connectedness shows that $\widetilde{f}(r)$
is independent of the path from $p$ to $r$ and we easily check that $%
\widetilde{f}$ is 1-1 and onto. Thus we have

\begin{proposition}
If $\widetilde{\varepsilon \mathcal{G}_{k}}$ is complete and simply
connected then it is globalizable.
\end{proposition}

If $\widetilde{\varepsilon \mathcal{G}_{k}}$ is complete but not
globalizable, then we can pull back $\widetilde{\varepsilon \mathcal{G}_{k}}$
to the universal covering space $\pi :\mathcal{M\rightarrow }M.$ Since $\pi
^{\ast }\widetilde{\varepsilon \mathcal{G}_{k}}$ is complete and $\mathcal{M}
$ is simple connected, $\pi ^{\ast }\widetilde{\varepsilon \mathcal{G}_{k}}$
globalizes to a Lie group $G$ acting on $\mathcal{M}$ and we obtain the
Klein geometry $(G,\mathcal{M)}\simeq G/H=\mathcal{M}$. To summarize, we have

\begin{proposition}
Let $\varepsilon \mathcal{G}_{k}$ be locally solvable and complete. Then the
pseudogroup $\widetilde{\varepsilon \mathcal{G}_{k}}$ globalizes to a Lie
group $G$ on the universal covering space $\mathcal{M}$ so that $(G,\mathcal{%
M)}\simeq G/H=\mathcal{M}$ and $M=G/H\smallsetminus \Gamma $ for some
discrete subgroup $\Gamma \subset G$ which is isomorphic to the fundamental
group of $M.$ \ 
\end{proposition}

A pseudogroup arising from a locally solvable $phg$ as above is a finite
type Lie pseudogroup according to [14].

Observe that we defined completeness of $\varepsilon \mathcal{G}_{k}$ only
when it is locally solvable. We will turn back to this issue in Section 5.

Conversely we now start with a transitive and effective Klein geometry $%
(G,M)\simeq G/H=M.$ We assume that $G$ is connected and $M$ is simply
connected (so $H$ is also connected) for reasons which will be clear below.
Note that $(G,M)$ always lifts to a Klein geometry $(\widehat{G},$ $\mathcal{%
M})$ where $\mathcal{M}$ is the universal cover of $M.$

We fix some base point $e\in M$ and let $H_{e}=\left\{ g\in G\mid
g(e)=e\right\} .$ Recall that a coordinate system around $e$ identifies $%
\mathcal{U}_{k}^{e,e}$ with the jet group $G_{k}(n).$ We have the evaluation
maps

\begin{eqnarray}
j_{i}^{e} &:&H_{e}\longrightarrow \mathcal{U}_{i}^{e,e}\simeq G_{i}(n) \\
&:&h\longrightarrow j_{i}(h)^{e,e}\text{ \ \ \ }0\leq i  \notag
\end{eqnarray}%
which are clearly homomorphisms of Lie groups. Since $G$ is connected and $%
(G,M)$ is effective (as we always assume in this note), there exists an
integer $k$ such that $j_{k}^{e}$ becomes injective ([2]).

\begin{definition}
The smallest integer $k$ such that (3) becomes injective is the order of $%
(G,M)$ denoted by $ord(G,M).$
\end{definition}

Since $G$ acts transitively, this integer is independent of our choice of
the base point $e.$ Clearly, $ord(G,M)=0$ if and only if $G$ acts simply
transitively.

Definition 8 needs only connectedness of $G.$ If $M$ is not simply
connected, then $ord(G,M)$ may be one greater than the one in Definition 8
and the length of the top filtration in (13) below may be one greater than
the bottom filtration.

Now any $g\in G$ is determined globally by its $k$-arrow $j_{k}(g)^{p,q}$
for any $p,q\in M.$ Indeed, $j_{k}(g)^{p,q}=j_{k}(g^{\prime })^{p,q}$ $%
\Leftrightarrow $ $j_{k}(g^{\prime }\circ g^{-1})^{p,p}$ $=Id_{k}^{p,p}$ $%
\Leftrightarrow g=g^{\prime }$ since $j_{k}^{p}$ is injective. We define the
groupoid $\mathcal{G}_{k}\subset \mathcal{U}_{k}$ on $M$ by defining its
fiber $\mathcal{G}_{k}^{p,q}\overset{def}{=}\left\{ j_{k}(f)^{p,q}\mid f\in
G,\text{ }f(p)=q\right\} .$ Further, $j_{k}(f)^{p,q}$ determines $f$ which
in turn determines $j_{k+1}(f)^{p,q}\overset{def}{=}\varepsilon
j_{k}(f)^{p,q}.$ Thus we obtain a splitting $\varepsilon $ such that $%
\mathcal{G}_{k}\simeq \varepsilon \mathcal{G}_{k}\subset \mathcal{U}_{k+1}.$
However there is a technical difficulty: Even though the map $j_{k}^{e}$ is
smooth as it is continuous, the image $j_{k}(H_{e})\subset $ $\mathcal{U}%
_{k}^{e,e}$ need not be a closed subgroup and therefore the groupoid $%
\varepsilon \mathcal{G}_{k}\subset \mathcal{U}_{k+1}$ need not be a
subgroupoid which should be an imbedded submanifold. Such an example is
given in [19]. If $H_{e}$ is compact, this anomaly can not occur but in this
case $ord(G,M)\leq 1$ by Proposition 23 below so this is a very strong
condition. We do not know any sufficient condition which makes $%
j_{k}(H_{e})\subset $ $\mathcal{U}_{k}^{e,e}$ closed but does not restrict $%
ord(G,M).$

In this note we will make the overall assumption

\textbf{A1:} The injection

\begin{equation}
j_{k}^{e}:H_{e}\rightarrow \mathcal{U}_{k}^{e,e}\simeq G_{k}(n)
\end{equation}%
imbeds $H_{e}$ as a closed subgroup for some (hence all) base point $e\in M$.

Clearly, $j_{k}(H_{e})\subset $ $\mathcal{U}_{k}^{e,e}$ is closed $%
\Leftrightarrow \varepsilon j_{k}(H_{e})\subset $ $\mathcal{U}_{k+1}^{e,e}$
is closed. Henceforth we will identify $H_{e}$ with its image $\varepsilon
j_{k}(H_{e}).$

Therefore we deduce

\begin{proposition}
A Klein geometry $(G,M)$ determines a locally solvable (in fact globally
solvable admitting $G$ as its global solution space) $phg$ $\varepsilon 
\mathcal{G}_{k}$ where $k=ord(G,M).$
\end{proposition}

To summarize what we have done so far, a locally solvable $\varepsilon 
\mathcal{G}_{k}$ makes $M$ locally homogeneous. With the assumption of
completeness, the universal cover $\widetilde{M}$ becomes globally
homogeneous. Conversely any Klein geometry $(G,M)$ with $ord(G,M)=k$
determines a globally solvable $\varepsilon \mathcal{G}_{k}$ with the above
assumptions$.$

Now suppose that the identification $\mathcal{U}_{k+1}^{e,e}\simeq
G_{k+1}(n) $ in (4) is induced by some coordinates $(U,x^{i})$ around $e.$
Since a change of coordinates $(x^{i})\rightarrow (y^{i})$ conjugates this
identification, $H_{e}\simeq \varepsilon j_{k}(H_{e})$ defines a unique
conjugacy class inside $G_{k+1}(n).$ Since $G$ acts transitively, this
conjugacy class is also independent of the choice of the basepoint $e.$ With
an abuse of notation, we denote this conjugacy class by $\left\langle
H\right\rangle _{G}$ where $H$ stands for any stabilizer of $(G,M)$.

\begin{definition}
The conjugacy class $\left\langle H\right\rangle _{G}$ inside $G_{k+1}(n)$
is the vertex class of $(G,M)\simeq G/H=M.$
\end{definition}

We now fix $H,$ $\dim M$ and want to understand the dependence of $%
\left\langle H\right\rangle _{G}$ from $G$ as $(G,M)$ varies where $G$ is
connected and $M$ is simply connected as we assumed above. The below
examples show that we may have $\left\langle H\right\rangle _{G}=$ $%
\left\langle H\right\rangle _{G^{\prime }}$ but $G,G^{\prime }$ are not even
locally isomorphic.

\textbf{Example 1: }Consider the projection $G_{1}(n)\rightarrow
G_{0}(n)=\{1\}$ with the only splitting $\varepsilon (1)=1.$ Let $%
\left\langle \varepsilon 1\right\rangle _{G_{1}(n)}=\{1\}$ denote the
conjugacy class of $\varepsilon G_{0}(n)$ inside $G_{1}(n).$ Now \textit{any 
}Klein geometry $(G,M)$ with $ord(G,M)=0$ defines this vertex class.

\textbf{Example 2: }Consider

\begin{equation}
0\rightarrow K_{2,1}\rightarrow G_{2}^{\circ }(n)\rightarrow G_{1}^{\circ
}(n)\rightarrow 1
\end{equation}%
where $G_{2}^{\circ }(n),$ $G_{1}^{\circ }(n)$ are the connected components
of $G_{2}(n),$ $G_{1}(n).$ Let $\varepsilon :G_{1}^{\circ }(n)\rightarrow
G_{2}^{\circ }(n)$ be the splitting defined by $(a_{j}^{i})\rightarrow
(a_{j}^{i},0)$, i.e., $a_{jk}^{i}=0,$ $1\leq i,j,k\leq n.$ Let $\left\langle
\varepsilon G_{1}^{\circ }(n)\right\rangle _{G_{2}(n)}$ denote conjugacy
class of $\varepsilon G_{1}^{\circ }(n)$ inside $G_{2}^{\circ }(n)$ which we
call the affine vertex class (of dimension $n$). Any other such splitting
defines the same conjugacy class! The Klein geometry $G_{1}^{\circ
}(n)\rtimes \mathbb{R}^{n}/G_{1}^{\circ }(n)$ has order one and

\begin{equation}
\left\langle G_{1}^{\circ }(n)\right\rangle _{G_{1}^{\circ }(n)\rtimes 
\mathbb{R}^{n}}=\left\langle \varepsilon G_{1}^{\circ }(n)\right\rangle
_{G_{2}^{\circ }(n)}
\end{equation}

\textbf{Example 3: }$\ $We restrict $\varepsilon $ in (5) to the orthogonal
group $SO(n)$ and let $\left\langle \varepsilon SO(n)\right\rangle
_{G_{2}(n)}$ denote the conjugacy class of $\varepsilon (SO(n))$ inside $%
G_{2}^{\circ }(n).$

Now consider the three Klein geometries $SO(n+1)/SO(n)=S^{n},$ $SO(n)\rtimes 
\mathbb{R}^{n}/SO(n)=\mathbb{R}^{n},$ $SO(n,1)/SO(n)=\mathbb{H}^{n}$. These
Klein geometries have order one and we have

\begin{equation}
\left\langle SO(n)\right\rangle _{SO(n+1)}=\left\langle SO(n)\right\rangle
_{SO(n)\rtimes \mathbb{R}^{n}}=\left\langle SO(n)\right\rangle
_{SO(n,1)}=\left\langle \varepsilon SO(n)\right\rangle _{G_{2}^{\circ }(n)}
\end{equation}

So three nonisomorphic Lie groups define the same vertex class.

\textbf{Example 4: }Consider

\begin{equation}
0\rightarrow K_{3,2}(1)\rightarrow G_{3}^{\circ }(1)\rightarrow G_{2}^{\circ
}(1)\rightarrow 1
\end{equation}

An element of $G_{3}(1)$ is an ordered triple $(a_{1},a_{2},a_{3})$, $%
a_{1}\neq 0$ and chain rule gives the group operation

\begin{equation}
(a_{1},a_{2},a_{3})(b_{1},b_{2},b_{3})=(a_{1}b_{1},a_{1}b_{2}+a_{2}(b_{1})^{2},a_{1}b_{3}+3a_{2}b_{1}b_{2}+a_{3}(b_{1})^{3})
\end{equation}%
We define $\varepsilon :G_{2}^{\circ }(1)\rightarrow G_{3}^{\circ }(1)$ by

\begin{equation}
\varepsilon (a_{1},a_{2})=(a_{1},a_{2},\frac{3(a_{2})^{2}}{2a_{1}})
\end{equation}%
Using (9) we check that $\varepsilon $ is a homomorphism (and $%
(a_{1},a_{2},\varepsilon (a_{1},a_{2}))^{-1}(a_{1},a_{2},a_{3})$ $%
=(1,0,S(a_{1},a_{2},a_{3}))$ where $S$ is the Schwarzian derivative!). Let $%
\left\langle \varepsilon G_{2}^{\circ }(1)\right\rangle _{G_{3}(1)}$ denote
the conjugacy class of $\varepsilon G_{2}^{\circ }(1)$ inside $G_{3}^{\circ
}(1).$

Let $\mathfrak{M}$ be the group of Mobius transformations $f(z)=\frac{az+b}{%
cz+d}$ normalized by $ad-bc=1$ acting transitively and effectively\textit{\ }%
on the sphere $S^{2}$. Now $ord(\mathfrak{M,}S^{2})=2$ and $\left\langle
H\right\rangle _{\mathfrak{M}}=\left\langle \varepsilon G_{2}^{\circ
}(1)\right\rangle _{G_{3}(1)}.$

The above examples show that the dependence of $\left\langle H\right\rangle
_{G}$ on $G$ is quite subtle. However the problem can be reduced to algebra
as follows. Let $(\mathfrak{g,h})$ be a pair of Lie algebras satisfying

$1)$ $\mathfrak{h\subset g}$ $\ $

$2)$ $\dim \mathfrak{g}-\dim \mathfrak{h}=n$

$3)$ $(\mathfrak{g,h)}$ is effective, i.e., $\mathfrak{h}$ contains no
nontrivial ideals inside $\mathfrak{g}.$

Now we fix $\mathfrak{h}$ and regard $\mathfrak{g}$ as variable. For any two
such pairs, we define $(\mathfrak{g,h)\sim }(\mathfrak{g}^{\prime }\mathfrak{%
,h)}$ if $\mathfrak{g\simeq g}^{\prime }$ and the isomorphism $\simeq $
restricts to identity on $\mathfrak{h}.$ The problem is to understand the
equivalence classes. In Example 1 this is the formidable problem of
classifying all Lie algebras whereas in Example 3 the solution is well known
from Riemannian geometry. By fixing $n$ and $\mathfrak{h}$, we call the
cardinality of the equivalence classes the uniformization number $\#(%
\mathfrak{h,}n).$ So $\#(0\mathfrak{,}n)=\infty ,$ $\#(\mathfrak{aff}(n)%
\mathfrak{,}n)=1$ and $\#(\mathfrak{o}(n)\mathfrak{,}n)=3,$ $n\geq 2.$ We
will comment more on $\#(\mathfrak{h,}n)$ in Appendix C.

We defined so far the vertex class $\left\langle H\right\rangle _{G}$ of $%
(G,M)\simeq G/H=M$. In the same way, we define the vertex class $%
\left\langle \varepsilon \mathcal{G}_{k}\right\rangle _{G_{k+1}(n)}$ of any $%
phg$ $\varepsilon \mathcal{G}_{k}$ as the conjugacy class of $\varepsilon 
\mathcal{G}_{k+1}^{p,p}\subset \mathcal{U}_{k+1}^{p,p}\simeq G_{k+1}(n)$
inside $G_{k+1}(n).$ Like $\left\langle H\right\rangle _{G},$ $\left\langle
\varepsilon \mathcal{G}_{k}\right\rangle _{G_{k+1}(n)}$ does not depend on
the identification $\mathcal{U}_{k+1}^{p,p}\simeq G_{k+1}(n)$ induced by
some coordinates around $p$ and is also independent of the choice of $p$ by
transitivity$.$

In the above examples we started with the conjugacy class of some subgroup $%
\varepsilon H$ inside some jet group and exhibited some Klein geometries ($%
\infty ,$ 1 and 3 in number) with the vertex classes equal to the conjugacy
class of $\varepsilon H$. Is this possible for any such $\varepsilon H?$ So
we face the following question

\textbf{Q: }For some arbitrary $phg$, does there exist some $G/H$ with $%
\left\langle \varepsilon \mathcal{G}_{k}\right\rangle
_{G_{k+1}(n)}=\left\langle H\right\rangle _{G}?$

We do not know the answer. Therefore we add the third condition $iii)$ to
the Definition 1:

$iii)$ There exists a Klein geometry $G/H$ with $\left\langle \varepsilon 
\mathcal{G}_{k}\right\rangle _{G_{k+1}(n)}=\left\langle H\right\rangle _{G}.$

$iii)$ is important for the following reason. We can restrict $\varepsilon 
\mathcal{G}_{k}$ to any open subset $U\subset M$ and the restriction $%
\varepsilon \mathcal{G}_{k\mid U}$ also satisfies $i),$ $ii).$ In Sections
6,7 we will assign certain invariants to $phg$'s which depend on their
equivalence class and vanish if there is a locally solvable $phg$ in this
equivalence class. We want these invariants vanish for $\varepsilon \mathcal{%
G}_{k\mid U}$ for any $\varepsilon \mathcal{G}_{k}$ if $U$ is sufficiently
small. This will be the case if $\varepsilon \mathcal{G}_{k\mid U}$ is
equivalent to some $phg$ on $U$ which is locally solvable and $iii)$ implies
this. In short, we want a $phg$ to be locally equivalent to a locally
solvable one.

According to $iii)$, we now state

\begin{equation}
\text{A }phg\text{ is modeled on some }\left\langle H\right\rangle _{G}\text{
defined by some }(G,M)\simeq G/H
\end{equation}

We now give some examples of $phg$'s on some (not necssarily simply
connected) $M.$

\textbf{Example 1 }(continues): For $k=0,$ $\mathcal{G}_{0}=M\times M$ and $%
\varepsilon $ assigns to any pair $(p,q)$ a unique $1$-arrow from $p$ to $q$%
. So $M$ admits $\varepsilon \mathcal{G}_{0}$ if and only if it is
parallelizable. Clearly, $\left\langle \varepsilon \mathcal{G}%
_{0}\right\rangle =\left\langle 1\right\rangle _{G}$ for any $G$ with $\dim
G=\dim M.$ If $\varepsilon \mathcal{G}_{0}$ is locally solvable, we get the
pseudogroup $\widetilde{\varepsilon \mathcal{G}_{0}}$ on $M$ which acts
simply transitively and, assuming completeness, globalizes to some Lie group 
$G$ on the universal cover of $M.$ Any Lie group $G$ is a possibility.
Therefore, the case $k=0$ gives the theory of parallelizable manifolds and
Lie groups as simply transitive transformation groups. In Section 10 we will
take a more careful look at this case.

\textbf{Example 2 }(continues): We recall $\mathcal{U}_{1}$ and let $%
\varepsilon $ be any symmetric connection on $T\rightarrow M.$ The
transformation rule of the components $\left( \varepsilon _{jk}^{i}\right) $
shows that $\varepsilon $ defines above any $1$-arrow of $\mathcal{U}_{1}$ a
unique $2$-arrow of $\mathcal{U}_{2}$. The $phg$ $\varepsilon \mathcal{U}%
_{1} $ has order one and $\left\langle \varepsilon \mathcal{U}%
_{1}\right\rangle =\left\langle \varepsilon G_{1}^{\circ }(n)\right\rangle
_{G_{2}(n)}$.

\begin{definition}
$\varepsilon \mathcal{U}_{1}$ is an affine geometry on $M.$
\end{definition}

The curvature $\mathcal{R}_{1}$ of $\varepsilon \mathcal{U}_{1},$ as defined
in Section 3, is \textit{not the same object }as the well known curvature $%
\mathcal{R}$ of $\varepsilon !$ Indeed $\mathcal{R}$ is a tensor whereas $%
\mathcal{R}_{1}$ is a second order object! However $\mathcal{R}_{1}=0$ $%
\Leftrightarrow \mathcal{R}=0.$ In this case, assuming completenes, the
pseudogroup $\widetilde{\varepsilon \mathcal{U}_{1}}$ globalizes to $%
Aff^{\circ }(\mathbb{R}^{n})=G_{1}^{\circ }(n)\rtimes \mathbb{R}^{n}$ on the
universal cover $\mathcal{M}\simeq \mathbb{R}^{n}$ of $M$. There is no
possibility other than $(Aff^{\circ }(\mathbb{R}^{n}),\mathbb{R}^{n}).$

\textbf{Example 3 }(continues): Let $g$ be a metric on $M.$ We define $%
\mathcal{G}_{1}^{p,q}$ as the set of all $1$-arrows from $p$ to $q$ which
map $g(p)$ to $g(q).$ The transformation rule of the Christoffel symbols $%
\varepsilon _{jk}^{i}$ shows that above any such $1$-arrow, there is a
unique $2$-arrow which defines $\varepsilon \mathcal{G}_{1}\subset \mathcal{U%
}_{2}.$ To be consistent with above general philosophy, we also assume that
the elements in the vertex groups have positive determinant.

\begin{definition}
$\varepsilon \mathcal{G}_{1}$ is a "Riemann geometry" on $M$
\end{definition}

Observe that a "Riemann geometry" according to Definition 12 is a second
order structure but not a first order structure! Now $\varepsilon \mathcal{G}%
_{1}$ is locally solvable $\Leftrightarrow $ $g$ has constant curvature. In
particular $\mathcal{R}_{1}$ is \textit{not the Riemann curvature tensor but
a second order geometric object}!! A. Blaom gives a very simple and explicit
formula for $\mathfrak{R}_{1}$ on pg. 6 of [4]. If $\varepsilon \mathcal{G}%
_{1}$ is locally solvable, the pseudogroup $\widetilde{\varepsilon \mathcal{G%
}_{1}}$ globalizes, assuming completeness, on the universal cover of $M$ to
one of the three groups in (7). There are no other possibilities other than
these three groups.

\textbf{Example 4} (continues)\textbf{: }As a generalization of Example 4,
we now want to define a projective geometry.

We first observe $\mathfrak{M\simeq }SL(2,\mathbb{R})$ $H\simeq $ the upper
triangular matrices = the stabilizer of $\infty $ obtained by setting $c=0.$

We will denote $H$ by $B(2).$ Now now fix $B(n)\subset SL(n,$ $\mathbb{R)}$
as our stabilizer. However we are not forced to fix $SL(n,$ $\mathbb{R)}$
because a $phg$ is not modeled on some $G/H$ but modeled only on $%
\left\langle H\right\rangle _{G}$ according to (11) (see the survey [10] for
the standard approach). For instance, we fix some entry just below the
diagonal and define $B(n)\subset P\subset SL(n,\mathbb{R)}$ by allowing only
that entry be nonzero below the diagonal. Then $P$ is a subgroup and $%
ord(P/B)=n.$ We can allow more than one entry below the diagonal but the
locations of these entries are not arbitrary. For $P=SL(n,\mathbb{R)}$,
however, $ord(P/B)=2$ (see (26) in [2]).

\begin{definition}
A projective geometry on $M$ is a $phg$ $\varepsilon \mathcal{G}_{k}$ on $M$
with vertex class $\left\langle B(n)\right\rangle _{P}$ for some Lie group $%
B(n)\subset P\subset SL(n,\mathbb{R)}$ as defined above.
\end{definition}

Observe that $k$ is determined by $\left\langle B(n)\right\rangle _{P}$ and $%
\dim M=\dim P-\dim M.$

In fact we can choose $P=G$ any Lie group $B(n)\subset G$ not necessarily
contained in $SL(n,\mathbb{R)}.$ It is our decision whether such a geometry
(if it exists at all as we require effectiveness) will qualify as a
"projective geometry". In the same way we define also a conformal geometry.
As long as we recognize the stabilizer, we can define \textit{that }geometry
with the freedom of choosing the vertex class. This process is very similar
to the classification of the principal bundles with some structure group $H$
and base $M$ and the decision of the vertex class may be interpreted as a
"geometrization" condition for the total space of that principal bundle
which is essentially a topological object.

We now turn back again to the imbedding (4). The filtration on the RHS of
(4) in terms of the projection of jets induces a filtration inside $G.$ Can
we recover this filtration group theoretically?

For some $(G,M)\simeq G/H$ we set $H=H_{0},$ $\mathfrak{g}=$ the Lie algebra
of $G$ and define inductively

\begin{equation}
H_{i+1}\overset{def}{=}\left\{ h\in H_{i}\mid Ad(h)x-x\in \mathfrak{h}_{i}%
\text{ \ for all }x\in \mathfrak{g}\text{, }i\geq 0\right\}
\end{equation}%
where $\mathfrak{h}_{i}$ denotes the Lie algebra of $H_{i}.$ Now $%
H_{i+1}\lhd H_{i}$ is a normal subgroup and we obtain the filtrations

\begin{eqnarray}
.... &\subset &H_{1}\subset H_{0}\subset G \\
... &\subset &\mathfrak{h}_{1}\subset \mathfrak{h}_{0}\subset \mathfrak{g} 
\notag
\end{eqnarray}

We can define the second filtration in (13) also as

\begin{equation}
\mathfrak{h}_{i+1}\overset{def}{=}\left\{ h\in \mathfrak{h}_{i}\mid \lbrack
h,x]\in \mathfrak{h}_{i}\text{ for all }x\in \mathfrak{g}\right\}
\end{equation}

If $0\neq \mathfrak{h}_{r}=\mathfrak{h}_{r+1}$ for some $r,$ then $\left\{
1\right\} \neq \cap _{i\geq r}H_{i}\lhd G$ which contradicts the
effectiveness of $G/H.$ The smallest integer $k$ such that $\mathfrak{h}%
_{k}=0$ is called the infinitesimal order of $G/H$ in [2]. Since $[\mathfrak{%
h}_{i},\mathfrak{h}_{j}]\subset \mathfrak{h}_{i+j}$ for $i,j\geq 0,$ $%
\mathfrak{h}_{i}/\mathfrak{h}_{i+1}$ is abelian for $i\geq 1$ and it is easy
to see that $\mathfrak{h}_{1}\subset Nil(\mathfrak{h})$ = the maximal
nilpotent ideal of $\mathfrak{h}$. In particular, $k=1$ if $\mathfrak{h}$ is
semisimple. If $H_{k}\neq \{1\}$ then $H_{k+1}=\{1\}$ because $%
H_{k+1}\subset Ker(Ad)=Z(G)$ $=$ the center of $G$ since $G$ is connected
which implies $H_{k+1}\subset Z(G)\cap H=\{1\}$ since $G/H$ is effective.
However, if $G/H=M$ is also simply connected then $H_{k}=\{1\}.$ Now using
(14) we define $ord(\mathfrak{g,h)}$ and obtain $ord(\mathfrak{g,h)=}%
ord(G/H) $ if $G/H=M$ is simply connected (which we assume henceforth in
this section). Now we have the commutative diagram

\begin{equation}
\begin{array}{ccccccc}
1\rightarrow & K_{r,s}(n) & \rightarrow & G_{r}(n) & \overset{\pi _{i,j}}{%
\rightarrow } & G_{s}(n) & \rightarrow 1 \\ 
& \uparrow &  & \uparrow j_{r} &  & \uparrow j_{s} &  \\ 
1\rightarrow & H_{s}/H_{r} & \rightarrow & H_{0}/H_{r} & \overset{\pi }{%
\rightarrow } & H_{0}/H_{s} & \rightarrow 1%
\end{array}%
\end{equation}%
where $r\leq s+1$ and the vertical imbeddings are induced by $j_{k}.$ The
principal bundle $G/H_{i}\rightarrow G/H_{j}$ with structure group $%
H_{j}/H_{i}$ can now be identified with the principal bundle $\mathcal{G}%
_{i}^{e,\bullet }\rightarrow \mathcal{G}_{j}^{e,\bullet }$ as

\begin{equation}
\begin{array}{ccc}
\mathcal{G}_{i}^{e,\bullet } & \rightarrow & \mathcal{G}_{j}^{e,\bullet } \\ 
\parallel &  & \parallel \\ 
G/H_{i} & \rightarrow & G/H_{j}%
\end{array}%
\end{equation}

From the definition of the filtration (14) we deduce

\begin{proposition}
$ord(\mathfrak{g,h}_{i})=k-i,$ $0\leq i\leq k$.
\end{proposition}

Proposition 14 together with (16) will play a fundamental role in Section 7.

\section{ The algebroid of a $phg$}

Since $\varepsilon \mathcal{G}_{k}$ is a groupoid, we can define its
algebroid by linearization. Our purpose is to describe this linearization
process in some detail. Since $\varepsilon \mathcal{G}_{k}\subset \mathcal{U}%
_{k+1}$ is a subgroupoid, this inclusion will hold also for the algebroids.
Therefore we will first recall the algebroid of $\mathcal{U}_{k+1}$ and
refer to [16], [17] for more details.

Let $T\rightarrow M$ be the tangent bundle and $J_{k}T\rightarrow M$ the $k$%
'th jet extension of $T\rightarrow M.$ Now $J_{k}T\rightarrow M$ is a vector
bundle whose fiber $\left( J_{k}T\right) ^{p}$ above $p\in M$ consists of
the $k$-jets of vector fields defined near $p.$ We have $J_{0}T=T.$ The
sections of $J_{k}T\rightarrow M$ are endowed with a bracket $[$ $,$ $]$
called the Spencer bracket which makes $J_{k}T\rightarrow M$ the algebroid
of $\mathcal{U}_{k}.$ To define $[$ $,$ $]$ we need two concepts. The first
is the ordinary Spencer operator

\begin{equation}
D:J_{k}T\rightarrow \wedge ^{\ast }T\otimes J_{k-1}T\text{ \ \ \ }k\geq 1
\end{equation}%
defined locally by the formula

\begin{equation}
(\xi ^{i},\xi _{j_{1}}^{i},...,\xi _{j_{k}...j_{1}}^{i})\rightarrow \left( 
\frac{\partial \xi ^{i}}{\partial x^{r}}-\xi _{r}^{i},...,\frac{\partial \xi
_{j_{k-1}...j_{1}}^{i}}{\partial x^{r}}-\xi _{rj_{k-1}...j_{1}}^{i}\right)
\end{equation}

The second is the algebraic bracket

\begin{equation}
\{\text{ },\text{ }\}_{p}:\left( J_{k}T\right) ^{p}\times \left(
J_{k}T\right) ^{p}\rightarrow \left( J_{k-1}T\right) ^{p}\text{ \ \ \ \ }%
k\geq 1
\end{equation}%
which is defined locally by differentiating the usual bracket formula $[\xi
(x),\eta (x)]^{i}(x)=$ $\xi ^{a}(x)\frac{\partial \eta ^{i}(x)}{\partial
x^{a}}-\eta ^{a}(x)\frac{\partial \xi ^{i}(x)}{\partial x^{a}}$ of two
vector fields $\xi =(\xi ^{i}(x))$, $\eta =(\eta ^{i}(x))$ $k$-times,
evaluating at $x=p$ and replacing all derivatives by jet variables. This
bracket does \textit{not} endow $J_{k}(T)^{p}$ with a Lie algebra structure
as it reduces the order of jets by one. However, let $\overline{J}%
_{k}T\rightarrow M$ denote the kernel of $J_{k}T\rightarrow J_{0}T$, that
is, the fiber $\left( \overline{J}_{k}T\right) ^{p}$ consists of all points
in $\left( J_{k}T\right) ^{p}$ which project to zero on the tangent space.
Now $\{$ $,$ $\}_{p}$ restricts to $\left( \overline{J}_{k}T\right)
^{p}\times \left( \overline{J}_{k}T\right) ^{p}\rightarrow \left( \overline{J%
}_{k}T\right) ^{p}$ and $\left( \overline{J}_{k}T\right) ^{p}$ endowed with $%
\{$ $,$ $\}_{p}$ is a Lie algebra. In fact, this Lie algebra is the Lie
algebra of the vertex group $\mathcal{U}_{k}^{p,p}.$ Clearly $\{$ $,$ $%
\}_{p} $ extends to a bracket on the sections of $J_{k}T\rightarrow M$ by
pointwise evaluation. We denote this bracket by $\{$ $,$ $\}.$

The Spencer bracket is defined now as follows. Let $\xi _{k},$ $\eta _{k}$
be two sections of $J_{k}T\rightarrow T.$ We lift $\xi _{k},$ $\eta _{k}$ to
some sections $\xi _{k+1},$ $\eta _{k+1}$ of $J_{k+1}T\rightarrow M$ and
define

\begin{equation}
\lbrack \xi _{k},\eta _{k}]\overset{def}{=}\{\xi _{k+1},\eta _{k+1}\}+i(\xi
_{0})D\eta _{k+1}-i(\eta _{0})D\xi _{k+1}\qquad k\geq 0
\end{equation}%
where $\xi _{0}$ is the projection of $\xi _{k}$ on the tangent space and $%
i(\xi _{0})$ denotes contraction with respect to $\xi _{0}.$ The Spencer
bracket $[\xi _{k},\eta _{k}]$ does not depend on the lifts. This backet
satisfies the Jacobi identity. Further it commutes with the projections $%
J_{k}T\rightarrow J_{r}T,$ $r\leq k$ and gives the usual bracket of vector
fields for $J_{0}T\rightarrow M.$ It also commutes with the prolongation of
vector fields: Some $\xi =(\xi ^{i})\in $ $\mathfrak{X}(M)$ prolongs to a
section of $J_{k}T\rightarrow T$ as $pr_{k}:(\xi ^{i})\rightarrow (\xi ^{i},%
\frac{\partial \xi ^{i}}{\partial x^{j_{1}}},...,\frac{\partial ^{k}\xi ^{i}%
}{\partial x^{j_{k}...j_{1}}})$ and we have $pr_{k}[\xi ,\eta ]=[pr_{k}\xi
,pr_{k}\eta ]$.

The above definition of the Spencer bracket is technical. To understand the
geometry behind it, we first observe that the vector bundle $%
J_{k}T\rightarrow M$ is associated with the groupoid $\mathcal{U}_{k+1}$ in
the sense that any $(k+1)$-arrow $j_{k+1}(f)^{p,q}$ induces an isomorphism

\begin{equation}
j_{k+1}(f)_{\ast }^{p,q}:\left( J_{k}T\right) ^{p}\rightarrow \left(
J_{k}T\right) ^{q}
\end{equation}%
defined locally by differentiationg the transformation rule $\frac{\partial
y^{i}}{\partial x^{a}}\eta ^{a}(x)=\eta ^{i}(y)$ of the vector field $\eta
=(\eta ^{i})$ $k$-times, evaluating at $x=p,$ $y=q,$ and replacing all
derivatives with jet variables. In particular (21) gives a faithful
representation of $\mathcal{U}_{k+1}^{p,p}\simeq G_{k+1}(n)$ on $\left(
J_{k}T\right) ^{p}.$ This representation descends to a representation of $%
\mathcal{U}_{k}^{p,p}\simeq G_{k}(n)$ on $\left( \overline{J}_{k}T\right)
^{p}\simeq \mathfrak{g}_{k}(n)$ which is the adjoint representation of $%
G_{k}(n)$ on its Lie algebra $\mathfrak{g}_{k}(n).$

Now we recall that the points of the principal bundle $\pi :\mathcal{U}%
_{k}^{e,\bullet }\rightarrow M$ are $k$-arrows eminating from the base point 
$e.$ So let $\overline{p}\in \mathcal{U}_{k}^{e,\bullet }$, $\pi (\overline{p%
})=p=j_{k}(f)^{e,p}$. Let $\xi _{\overline{p}}$ be a tangent vector at $%
\overline{p}$ which projects to the tangent vector $\xi _{p}$ at $p.$ By
acting with the structure group $\mathcal{U}_{k}^{e,e}$ on the fiber $%
\mathcal{U}_{k}^{e,p}$ we translate $\xi _{\overline{p}}$ to the points in
the fiber $\pi ^{-1}(p)=\mathcal{U}_{k}^{e,p}.$ We call this data a set of
parallel vectors at the fiber $\pi ^{-1}(p).$ We now have

\begin{proposition}
There is a canonical identification between the following objects.

$i)$ The set of parallel vectors at the fiber $\pi ^{-1}(p)$

$ii)$ The fiber $\left( J_{k}T\right) ^{p}$ of the vector bundle $%
J_{k}T\rightarrow M$ over $p.$
\end{proposition}

In particular, Proposition 15 gives the canonical identification

\begin{equation}
T(\mathcal{U}_{k}^{e,\bullet })^{\overline{p}}\simeq \left( J_{k}T\right)
^{p}
\end{equation}%
where $T(\mathcal{U}_{k}^{e,\bullet })^{\overline{p}}$ denotes the tangent
space of $\mathcal{U}_{k}^{e,\bullet }$ at $\overline{p}.$

The reason for (22) is simple: Let $\overline{p}\in \mathcal{U}%
_{k}^{e,\bullet }$, $\pi (\overline{p})=p=j_{k}(f)^{e,p}$. A diffeomorphism $%
g$ on $M$ lifts to a diffeomorphism $\overline{g}$ on $\mathcal{U}%
_{k}^{e,\bullet }$ defined by $j_{k}(f)^{e,p}\rightarrow j_{k}(g\circ
f)^{e,g(p)}.$ Consider the $1$-parameter group $g_{t}(x)$ of local
diffeomorphisms defined by a vector field $\xi (x)$ defined around $p.$ Now $%
\xi (x)$ lifts to a vector field on $\mathcal{U}_{k}^{e,\bullet }$ whose
value at $\overline{p}$ depends on $j_{k}(\xi )^{p}$ and this map is an
isomorphism.

Proposition 15 gives the canonical identification

\begin{equation}
\left\{ \text{right invariant vector fields on }\mathcal{U}_{k}^{e,\bullet
}\rightarrow M\right\} \simeq \left\{ \text{sections of }J_{k}T\rightarrow
M\right\}
\end{equation}

Since the LHS of (23) is a Lie algebra with the usual bracket of vector
fields on $\mathcal{U}_{k}^{e,\bullet },$ we get a bracket on the RHS of
(23)...which is the Spencer bracket.

The linearization of the nonlinear $PDE$ $\varepsilon \mathcal{G}_{k}$ to
the linear $PDE$ $\varepsilon \mathfrak{G}_{k}$ is best understood in
coordinates. We replace the "finite transformations" in (2) by
"infinitesimal transformations", i.e., we substitute $y^{i}=x^{i}+t\xi ^{i},$
$f_{j_{1}}^{i}=\delta _{j_{1}}^{i}+t\xi _{j_{1}}^{i},$ $f_{j_{2}j_{1}}^{i}=t%
\xi _{j_{2}j_{1}}^{i},$ ..., $f_{j_{k+1}...j_{1}}^{i}=t\xi
_{j_{k+1}...j_{1}}^{i}$ into (2) and differentiate at $t=0.$ The resulting
equations

\begin{equation}
\widehat{\Phi ^{\alpha }}:(x^{i},\xi ^{i},...,\xi _{j_{k+1}...j_{1}}^{i})=0
\end{equation}%
are linear in the variables $\xi ^{i},...,\xi _{j_{k+1}...j_{1}}^{i}$ which
are the local coordinates on $J_{k+1}T$ over the fiber $\pi ^{-1}(x).$ As in
(2), the top coordinates $\xi _{j_{k+1}...j_{1}}^{i}$ can be solved uniquely
in terms of $\xi ^{i},...,\xi _{j_{k}...j_{1}}^{i}$ (we will denote this
splitting again by $\varepsilon )$ and (24) puts no restriction on the
variables $x^{i}.$ So (24) defines a subbundle $\varepsilon \mathfrak{G}%
_{k}\rightarrow M$ of $J_{k+1}T\rightarrow M.$ The crucial fact is that the
Spencer bracket restricts to the sections of $\varepsilon \mathfrak{G}%
_{k}\rightarrow M.$ Since $\mathcal{G}_{k}\rightarrow $ $\varepsilon 
\mathcal{G}_{k}$ is an isomorphism of groupoids, $\mathfrak{G}%
_{k}\rightarrow \varepsilon \mathfrak{G}_{k}$ is an isomorphism of
algebroids, i.e., it preserves the bracket. Thus we have the diagram

\begin{equation}
\begin{array}{ccc}
\mathcal{U}_{k+1} & \Longrightarrow & J_{k+1}T\rightarrow M \\ 
\cup &  & \cup \\ 
\varepsilon \mathcal{G}_{k} & \Longrightarrow & \varepsilon \mathfrak{G}%
_{k}\rightarrow M%
\end{array}%
\end{equation}%
where $\Longrightarrow $ denotes linearization.

At this point, it is possible to define the Lie algebroid $\varepsilon 
\mathfrak{G}_{k}\rightarrow M$ independently as an "infinitesimal $phg$ of
order $k$ on $M$". However, once properly defined, this object will be the
linearization of a unique "finite $phg$ $\varepsilon \mathcal{G}_{k}$".

Like the Spencer bracket, all the calculus on $J_{k+1}T\rightarrow M$ (which
we only touched here) restricts to $\varepsilon \mathfrak{G}_{k}\rightarrow
M.$ For instance, for $j_{k+1}(f)^{p,q}\in \varepsilon \mathcal{G}%
_{k}^{p,q}, $ (21) restricts as

\begin{equation}
\varepsilon j_{k+1}(f)_{\ast }^{p,q}:\mathfrak{G}_{k}^{p}\rightarrow 
\mathfrak{G}_{k}^{q}
\end{equation}%
where $\mathfrak{G}_{k}^{p}$ denotes the fiber of $\mathfrak{G}%
_{k}\rightarrow M$ over $p.$ Therefore, even though $J_{k}T\rightarrow M$ is
associated with $\mathcal{U}_{k+1}$, $\mathfrak{G}_{k}\rightarrow M$ is
associated with $\mathcal{G}_{k}.$ This is a particular case of the
"stabilization of the order of jets using the splitting $\varepsilon $"
which will play a fundamental role in this note. Similarly, the
identification (22) restricts as

\begin{equation}
T(\mathcal{G}_{k}^{e,\bullet })^{\overline{p}}\simeq \mathfrak{G}_{k}^{p}
\end{equation}

Now (24) shows that the points $\overline{p}$ in the fiber $\mathfrak{G}%
_{k}^{p}$ of $\varepsilon \mathfrak{G}_{k}\rightarrow M$ over $p$ are
"initial conditions" for the linear $PDE$ $\varepsilon \mathfrak{G}%
_{k}\rightarrow M$ locally defined by (24). We call $\varepsilon \mathfrak{G}%
_{k}\rightarrow M$ locally solvable at $p\in M$ if for any $\overline{p}\in 
\mathfrak{G}_{k}^{p}$ there exists a vector field $\xi $ defined near $p$
whose prolongation $pr_{k+1}(\xi )$ is a section of $\varepsilon \mathfrak{G}%
_{k}\rightarrow M$ passing through $\overline{p}.$

\begin{definition}
$\varepsilon \mathfrak{G}_{k}\rightarrow M$ is locally solvable if it is
locally solvable on $M.$
\end{definition}

If $\varepsilon \mathfrak{G}_{k}\rightarrow M$ is locally solvable, then its
solutions are determined locally by their initial conditions. Therefore
"analytic continuation" is possible along paths. Note that the Spencer
bracket becomes the ordinary bracket of vector fields on local solutions.
Thus we can define the presheaf of Lie algebras $\mathfrak{g}(U)\overset{def}%
{=}$ the Lie algebra of local solutions on $U.$

Now we have the following fundamental

\begin{proposition}
$\varepsilon \mathcal{G}_{k}\rightarrow M$ is locally solvable $%
\Longleftrightarrow \varepsilon \mathfrak{G}_{k}\rightarrow M$ is locally
solvable.
\end{proposition}

The implication $\Rightarrow $ follows easily from definitions whereas $%
\Leftarrow $ is quite nontrivial. To see what is involved in Proposition 17,
assume that $\varepsilon \mathcal{G}_{k}\rightarrow M$ is locally solvable
and the pseudogroup $\widetilde{\varepsilon \mathcal{G}_{k}}$ globalizes to $%
G$ so that $\mathcal{G}_{k}^{e,\bullet }\simeq G.$ This implies that the
local solutions $\widetilde{\varepsilon \mathfrak{G}_{k}}$ of $\varepsilon 
\mathfrak{G}_{k}\rightarrow M$ also globalize and we obtain a Lie algebra $%
\mathfrak{g}$ of vector fields on $M.$ Not surprisingly, $\mathfrak{g}$ is
the Lie algebra of the infinitesimal generators of the Klein geometry $%
(G,M). $ Since $(G,M)$ is effective by construction, $\mathfrak{g}$ is
isomorphic to the Lie algebra of the abstract Lie group $G.$

So the bottom line of (25) becomes the assignment

\begin{equation}
\left\{ \text{The Klein geometry }(G,M)\right\} \Rightarrow \left\{ \text{%
the Lie algebra }\mathfrak{g}\text{ of the infinitesimal generators}\right\}
\end{equation}

So $\Leftarrow $ of Proposition 17 asserts that locally $(G,M)$ can
recovered from $\mathfrak{g.}$ For an abstract Lie group $G$, observe that
the assignment $\left\{ G\right\} \Rightarrow \left\{ \text{its Lie algebra }%
\mathfrak{g}\right\} $ involves a choice of left/right and is not canonical
whereas (28) is canonical even for $k=0,$ i.e., if $G$ acts simply
transitively. In this simplest case, $\Leftarrow $ of Proposition 17 becomes
the classical version (not the Cartan's version) of the Lie's 3rd Theorem.

Finally we note that, assuming local solvability, \textbf{A1 }linearizes to
an injection of Lie algebras and we get a correspondence between the
algebraic filtration in (14) and the "jet filtration" on $\mathfrak{G}%
_{k}\rightarrow M$ in terms of the projection of jets. In particular we get
the linearizations of (15), (16).

\section{Curvature}

Let $\varepsilon \mathcal{G}_{k}$ be a $phg$ of order $k$ on $M.$ In
Sections 2, 3 we have seen that the local solvability of $\varepsilon 
\mathcal{G}_{k}$ and its algebroid $\varepsilon \mathfrak{G}_{k}$ is a
fundamental concept because the theory of locally solvable $phg$'s is the
same as the theory of homogeneous spaces. Since local solvability is a very
intuitive concept, we can easily guess at some constructions and theorems
assuming local solvability without going into technical proofs. Since this
is a qualitative concept, it is desirable to define a quantity $\mathcal{R}%
_{k}$ = the curvature of $\varepsilon \mathcal{G}_{k}$ in such a way that we
will have

\begin{equation}
\mathcal{R}_{k}=0\Longleftrightarrow \varepsilon \mathcal{G}_{k}\text{ is
locally solvable}
\end{equation}

Similiarly we want to define $\mathfrak{R}_{k}=$ the curvature of $%
\varepsilon \mathfrak{G}_{k}$ such that

\begin{equation}
\mathfrak{R}_{k}=0\Longleftrightarrow \varepsilon \mathfrak{G}_{k}\text{ is
locally solvable}
\end{equation}

In view of Proposition 17 we will have

\begin{equation}
\mathcal{R}_{k}=0\Longleftrightarrow \mathfrak{R}_{k}=0
\end{equation}

Further, since $\varepsilon \mathfrak{G}_{k}$ is the linearization of $%
\varepsilon \mathcal{G}_{k}$, we require that $\mathfrak{R}_{k}$ should be
obtained from $\mathcal{R}_{k}$ by the same linearization process.

Probably the first thing that comes to mind is the following: Consider the
first prolongation $J_{1}\mathcal{G}_{k}^{e,\bullet }\rightarrow $ $\mathcal{%
G}_{k}^{e,\bullet }$ of the principal bundle $\mathcal{G}_{k}^{e,\bullet
}\rightarrow M.$ Sections of $J_{1}\mathcal{G}_{k}^{e,\bullet }\rightarrow 
\mathcal{G}_{k}^{e,\bullet }$ are in 1-1 correspondence with connections on $%
\mathcal{G}_{k}^{e,\bullet }\rightarrow M.$ Now it is easy to show

\begin{proposition}
The splitting $\mathcal{G}_{k}\rightarrow \varepsilon \mathcal{G}_{k}$
defines a connection on $\mathcal{G}_{k}^{e,\bullet }\rightarrow M.$
\end{proposition}

It is natural to define $\mathcal{R}_{k}$ to be the curvature of this
connection. It is an extremely surprising fact (probably more than that!)
that $\mathcal{R}_{k}$ is not this curvature!! There is a very short
conceptual way of seeing this as follows: For simplicity we assume $\mathcal{%
R}_{k}\mathcal{=}0$ and $\widetilde{\varepsilon \mathcal{G}_{k}}$ globalizes
to $G$ so that $\mathcal{G}_{k}^{e,\bullet }\simeq G$ and $G$ \textit{acts }%
transitively on $M.$ However in the general theory of principal bundles, the
total space of the principal bundle, in particular $\mathcal{G}%
_{k}^{e,\bullet }\simeq G$ in our case, does not act on the base
manifold...so these two curvatures can not be the same objects in general.
This will follow also from the technical definition (35) of $\mathcal{R}_{k}$
below. At this point it is crucial to observe that the principal bundle and
the connection are seperate objects in the general theory whereas they unify
into a single object in the definition of a $phg$. In particular, $\mathcal{R%
}_{k}$ is not the curvature of a connection on any principal bundle but is
the curvature of $\varepsilon \mathcal{G}_{k}.$

To find the technical definition of $\mathcal{R}_{k},$ we consider the first
nonlinear Spencer sequence ([12], [16], [17])

\begin{equation}
\mathbf{1\longrightarrow }Aut(M)\overset{j_{k}}{\longrightarrow }\mathcal{U}%
_{k}\overset{D_{1}}{\longrightarrow }T^{\ast }\otimes J_{k-1}T\overset{D_{2}}%
{\longrightarrow }\wedge ^{2}T^{\ast }\otimes J_{k-2}T
\end{equation}%
The explicit local formulas describing $D_{1},$ $D_{2}$ are given in [17],
pg.213-216 and it is quite easy to do computations with these formulas (if
we have enough patience!). Note that $k\geq 2$ in (32). Observe that $D_{1},$
$D_{2}$ reduce the order of jets by one.

Now (32) restricts as

\begin{equation}
\mathcal{G}_{k}\overset{D_{1}}{\longrightarrow }T^{\ast }\otimes \mathfrak{G}%
_{k-1}\overset{D_{2}}{\longrightarrow }\wedge ^{2}T^{\ast }\otimes \mathfrak{%
G}_{k-2}
\end{equation}%
and still $k\geq 2$ in (33). The crucial fact is that the splittings $%
\varepsilon :$ $\mathcal{G}_{k}\rightarrow \varepsilon \mathcal{G}_{k}$ and $%
\varepsilon :\mathfrak{G}_{k}\rightarrow \varepsilon \mathfrak{G}_{k}$
stabilize the order of jets in (33) as

\begin{equation}
\mathcal{G}_{k}\overset{D_{1}^{\prime }}{\longrightarrow }T^{\ast }\otimes 
\mathfrak{G}_{k}\overset{D_{2}^{\prime }}{\longrightarrow }\wedge
^{2}T^{\ast }\otimes \mathfrak{G}_{k}
\end{equation}%
and now $k\geq 0$ in (34). Even though $D_{2}\circ $ $D_{1}=0$ in (33), we
do not have $D_{2}^{\prime }\circ $ $D_{1}^{\prime }=0$ in (34). We define

\begin{equation}
\mathcal{R}_{k}\overset{def}{=}D_{2}^{\prime }\circ D_{1}^{\prime }
\end{equation}

So for any $k$-arrow $\alpha _{k}^{p,q}\in \mathcal{G}_{k}^{p,q}$ and $\xi
_{p},\eta _{p}\in T_{p}(M)$, we have

\begin{equation}
\mathcal{R}_{k}(\alpha _{k}^{p,q})(\xi _{p},\eta _{p})\in \mathfrak{G}%
_{k}^{q}
\end{equation}%
where, as before, $\mathfrak{G}_{k}^{q}$ denotes the fiber of the algebroid $%
\mathfrak{G}_{k}\rightarrow M$ over $q\in M.$

We now have

\begin{proposition}
The following are equivalent

$i)$ $\mathcal{R}_{k}\mathcal{=}0$

$ii)$ $\varepsilon \mathcal{G}_{k}$ is locally solvable

$iii)$ (34) is locally exact at $T^{\ast }\otimes \mathfrak{G}_{k}$

$iv)$ $\varepsilon \mathcal{G}_{k}$ is involutive
\end{proposition}

If one of the conditions of Proposition 19 holds, then (34) extends to the
second nonlinear Spencer sequence

\begin{equation}
\mathbf{1}\longrightarrow \widetilde{\varepsilon \mathcal{G}_{k}}\overset{%
j_{k}}{\longrightarrow }\mathcal{G}_{k}\overset{D_{1}^{\prime }}{%
\longrightarrow }T^{\ast }\otimes \mathfrak{G}_{k}\overset{D_{2}^{\prime }}{%
\longrightarrow }\wedge ^{2}T^{\ast }\otimes \mathfrak{G}_{k}
\end{equation}%
which is locally exact. It is instructive to check (37) in the simplest case 
$k=0$ of parallelizable manifolds studied in detail in [1] and construct
(37) in this case. Observe that there is no curvature in (37) (see however
[17], pg. 216 where the operator $D_{2}^{\prime }$ is claimed to be the
curvature).

Now we come to $\mathfrak{R}_{k}.$ Let $\overline{p}=\alpha _{k}^{e,p}$ and $%
\overline{q}=\beta _{k}^{e,q}$ be two points of $\mathcal{G}_{k}^{e,\bullet
}.$ According to (36) we have

\begin{equation}
\mathcal{R}_{k}(\overline{q}\circ \overline{p}^{-1})(\xi _{p},\eta _{p})\in 
\mathfrak{G}_{k}^{q}
\end{equation}

We now fix $\overline{p},\xi _{p},\eta _{p}$ in (38) and let $\overline{q}$
approach $\overline{p}$ along the direction of some tangent vector $\sigma
_{k}^{p}\in \mathfrak{G}_{k}^{p}=T(\mathcal{G}_{k}^{e,\bullet })^{\overline{p%
}}.$ The limiting value is an element of $\mathfrak{G}_{k}^{p}$ which we
write as

\begin{equation}
\left( \mathfrak{R}_{k}(\overline{p})(\xi _{p},\eta _{p})\right) (\sigma
_{k}^{p})
\end{equation}

The function $\sigma _{k}^{p}\rightarrow \left( \mathfrak{R}_{k}(\overline{p}%
)(\xi _{p},\eta _{p})\right) (\sigma _{k}^{p})$ turns out to be linear and
therefore

\begin{equation}
\mathfrak{R}_{k}(\overline{p})(\xi _{p},\eta _{p})\in Hom(\mathfrak{G}%
_{k}^{p},\mathfrak{G}_{k}^{p})
\end{equation}

So the object $\mathfrak{R}_{k}$ assigns to a point $\overline{p}$ on the
principal bundle $\mathcal{G}_{k}^{e,\bullet }\rightarrow M$ and two tangent
vectors $\xi _{p},\eta _{p}$ at $p\in M$ a linear map on the tangent space $%
T(\mathcal{G}_{k}^{e,\bullet })^{\overline{p}}$. Observe that $\mathfrak{R}%
_{k}$ is not an ordinary $2$-form on $M$.

Since (40) is the linearization of (36) which arises from the second
nonlinear Spencer sequence, we should be able to derive $\mathfrak{R}_{k}$
directly from the second linear Spencer sequence. Indeed, the ordinary
Spencer operator (17) restricts to

\begin{equation}
\mathfrak{G}_{k}\overset{D}{\longrightarrow }T^{\ast }\otimes \mathfrak{G}%
_{k-1}
\end{equation}%
where $k\geq 1.$ The splitting $\varepsilon :\mathfrak{G}_{k}\rightarrow
\varepsilon \mathfrak{G}_{k}$ stabilizes the order of jets in (41) as

\begin{equation}
\mathfrak{G}_{k}\overset{D^{\prime }}{\longrightarrow }T^{\ast }\otimes 
\mathfrak{G}_{k}
\end{equation}

Acting with $d$ on $T^{\ast }$ and with $D^{\prime }$ on $\mathfrak{G}_{k}$
we extend (42) one step the right as

\begin{equation}
\mathfrak{G}_{k}\overset{D^{\prime }}{\longrightarrow }T^{\ast }\otimes 
\mathfrak{G}_{k}\overset{D^{\prime \prime }}{\longrightarrow }\wedge
^{2}T^{\ast }\otimes \mathfrak{G}_{k}
\end{equation}%
where $D^{\prime \prime }=D_{2}^{\prime }$ in (37). Clearly we have

\begin{equation}
\mathfrak{R}_{k}=D^{\prime \prime }\circ D^{\prime }
\end{equation}

Now (43) extends to the sequence

\begin{equation}
\mathfrak{G}_{k}\overset{D^{\prime }}{\longrightarrow }T^{\ast }\otimes 
\mathfrak{G}_{k}\overset{D^{\prime \prime }}{\longrightarrow }\wedge
^{2}T^{\ast }\otimes \mathfrak{G}_{k}\longrightarrow ....\longrightarrow
\wedge ^{n}T^{\ast }\otimes \mathfrak{G}_{k}
\end{equation}%
which is locally exact if $\mathfrak{R}_{k}=0$ and (45) is the linear second
Spencer sequence. To summarize, $\mathcal{R}_{k}$ and $\mathfrak{R}_{k}$ are
obstructions to the passage from the first to the second Spencer sequences.

Finally, let $P=\mathcal{G}_{k}^{e,\bullet },$ $H=\mathcal{G}_{k}^{e,e}$ and 
$\mathfrak{h}$ the Lie algebra of $H.$ A connection on $P\rightarrow M$ is
an $\mathfrak{h}$-valued $1$-form on $P$ and its curvature $R$ is an $%
\mathfrak{h}$-valued $2$-form on $P$. Assume that $\mathfrak{R}_{k}=0$ and $%
\mathcal{G}_{k}^{e,\bullet }$ globalizes to $G$ which acts transitively on $%
M.$ Let $\mathfrak{g}$ be the Lie algebra of $G.$ Now (40) shows that $%
\mathfrak{R}_{k}$ is a $2$-form on $P$ with values in $Hom(\mathfrak{g},%
\mathfrak{g)}$ and therefore $\mathfrak{R}_{k}\neq R$ in general (however,
see Section 10 for the remarkable case $k=0!$)$\mathfrak{.}$ Note that $%
\mathfrak{R}_{k}$ can not be also the curvature of a Cartan connection on $%
\mathcal{G}_{k}^{e,\bullet }\rightarrow M$ which is a $\mathfrak{g}$-valued $%
1$-form on $\mathcal{G}_{k}^{e,\bullet }$ and its curvature is a $\mathfrak{g%
}$-valued $2$-form. However, we will see in Section 9 that \textit{this will
be the case }with a rather restrictive assumption.

\section{Cartan algebroids}

In [3], [4], Blaom proposed a very interesting and general theory of
infinitesimal geometric structures. The general philosophy, which he
attributes to E. Cartan, is to view an infinitesimal geometric structure as
a symmetry deformed by curvature. For this purpose, he defines the concept
of the Cartan algebroid. These are algebroids equipped with a linear
connection (which he calls somewhat confusingly the Cartan connection) whose
covariant derivative is compatible with the algebroid structure. The
curvature of the Cartan connection vanishes if and only if $M$ is locally
homogeneous. As an important fact, a Cartan algebroid is defined without the
use of jets and need not be transitive!

We now have

\begin{proposition}
The algebroid $\varepsilon \mathfrak{G}_{k}\rightarrow M$ of the $phg$ $%
\varepsilon \mathcal{G}_{k}$ is a transitive Cartan algebroid. The first
operator $D^{\prime }$ in (45) is the Cartan connection of $\varepsilon 
\mathfrak{G}_{k}\rightarrow M$ and $\mathfrak{R}_{k}$ defined by (44) is its
curvature.
\end{proposition}

We believe that all transitive Cartan algebroids arise as the Lie algebroids
of $phg$'s (possibly with some mild conditions of regularity). Therefore, we
believe that the theory of Cartan groupoids, whose study is initiated in the
Appendix A of\ [4] and expanded in [6], is essentially the same as the
theory of $phg$'s \textit{in the transitive case.}

In the recent preprint [5], Blaom also clarifies the the concept of
completeness of a not necessarily flat Cartan algebroid. We believe that
this will have important consequences for the theory of $PDE$'s in view of
the "equivalence" of transitive Cartan algebroids and $phg$'s.

\section{Characteristic classes on the base}

One of the great achievements of global the differential geometry in this
century is the theory of characteristic classes on principal and vector
bundles. These classes are cohomology classes on the base manifold which
measure the twisting of the bundle, i.e., their deviation from being
globally trivial. This is a topological theory. Therefore it came as a great
surprise when in 1970 R. Bott showed that Chern classes are also
obstructions to integrability of the to plane fields, i.e., subbundles of
the tangent bundle.

Now let $\varepsilon \mathcal{G}_{k}$ be a $phg$ of order $k$ and consider
the algebroid $\mathfrak{G}_{k}\rightarrow M.$ Let $\mathcal{P}^{\ast }(M,%
\mathfrak{G}_{k})$ denote the Pontryagin algebra ($P$-algebra for short) of
the vector bundle $\mathfrak{G}_{k}\rightarrow M.$ To recall the definition
of $\mathcal{P}^{\ast }(M,\mathfrak{G}_{k}),$ let $m$ denote the fiber
dimension of $\mathfrak{G}_{k}\rightarrow M.$ Let $\mathfrak{gl}(m,\mathbb{R)%
}$ denote the Lie algebra of $GL(m,\mathbb{R)}=G_{1}(m)=$ the Lie group of
invertible $m\times m$ matrices. So $\mathfrak{gl}(m,\mathbb{R)}$ is the
linear Lie algebra of all $m\times m$ matrices. A polynomial function $%
\varphi :\mathfrak{gl}(m,\mathbb{R)\rightarrow R}$ is called invariant if $%
\varphi (gAg^{-1})=\varphi (A)$ for all $A\in \mathfrak{gl}(m,\mathbb{R)},$ $%
g\in GL(m,\mathbb{R)}.$ The vector space $I_{GL(m,\mathbb{R)}}$ of all
invariant polynomials is an algebra generated by $T_{j}(A)\overset{def}{=}%
Trace(A^{j}),$ $j\geq 0.$ Now any connection $\omega $ on the vector bundle $%
\mathfrak{G}_{k}\rightarrow M$ defines the algebra homomorphism $CW$

\begin{eqnarray}
CW &:&I_{GL(m,\mathbb{R)}}\longrightarrow H_{dR}^{\ast }(M,\mathbb{R)} \\
&:&\varphi \longrightarrow \varphi (\kappa )  \notag
\end{eqnarray}%
where $\kappa $ is the curvature of $\omega .$ The map (46) is independent
of the connection$.$ If $\varphi $ is homogeneous of degree $r$, then $%
\varphi (\kappa )\in H_{dR}^{2r}(M,\mathbb{R)}$. Now $\mathcal{P}^{\ast }(M,%
\mathfrak{G}_{k})$ is the image of $CW$ and it can be shown (see [7]) that $%
\mathcal{P}^{j}(M,\mathfrak{G}_{k})=0$ if $j$ is not divisible by $4.$

As we observed in Section 3, the splitting $\mathfrak{G}_{k}\rightarrow
\varepsilon \mathfrak{G}_{k}$ defines a the particular linear connection $%
\varepsilon $ on the vector bundle $\mathfrak{G}_{k}\rightarrow M$ with
curvature $\mathfrak{R}_{k}\mathfrak{.}$ Further, $\mathfrak{R}%
_{k}=0\Leftrightarrow \mathcal{R}_{k}=0\Leftrightarrow $ $\varepsilon 
\mathcal{G}_{k}$ is locally solvable by Propositions 17,19. Since the map
(46) is independent of the connection, we obtain

\begin{proposition}
If $\varepsilon \mathcal{G}_{k}$ is locally solvable, then $\mathcal{P}%
^{\ast }(M,\mathfrak{G}_{k})=0.$
\end{proposition}

If $k=0$, note that $\mathcal{P}^{\ast }(M,T)=0$ without the assumption of
local solvability of $\varepsilon \mathcal{G}_{0}$ since the existence of $%
\varepsilon \mathcal{G}_{0}$ is equivalent to the parallelizability of $M.$

Proposition 21 follows almost trivially from our definitions. However it
gives a totally new way of looking at characteristic classes: among the set
of all vector bundles over $M,$ there is a particular subset of vector
bundles with the property that the restriction of the functor $\mathcal{P}%
^{\ast }$ to this subset gives global obstructions to integrability in the
sense of local solvability.

Recalling that $\mathfrak{G}_{0}=T,$ we have the exact sequence of vector
bundles

\begin{equation}
0\longrightarrow \mathcal{I}_{k}\longrightarrow \mathfrak{G}%
_{k}\longrightarrow T\longrightarrow 0
\end{equation}

For $k=0$, $\mathfrak{G}_{0}=T$ and $\mathcal{I}_{0}=0.$ For $k\geq 1,$ $%
\mathcal{I}_{k}\rightarrow M$ is a bundle of Lie algebras whose fiber over $%
p $ consists of all $k$-jets of vector fields at $p$ which project to zero
on the tangent space at $p$ and this fiber is the Lie algebra of the vertex
group $\mathcal{G}_{k}^{p,p}$. From (47) we conclude $\mathfrak{G}_{k}=%
\mathcal{I}_{k}\oplus T.$ Therefore if $\varepsilon \mathcal{G}_{k}$ is
locally solvable, then $p(\mathcal{I}_{k})\cdot p(T)=1$ which is the first
indication that the existence of some locally solvable $\varepsilon \mathcal{%
G}_{k}$ puts restrictions on $\mathcal{P}^{\ast }(M,T).$ To dig this point
deeper, we will recall some facts from [9] which have an intriguing relation
to Proposition 21.

Let $J_{(m)}T\rightarrow M$ denote the $m$-th order iterated jet bundle of $%
T\rightarrow M,$ i.e., $J_{(m)}T\rightarrow M$ is obtained from $%
T\rightarrow M$ by applying the functor $J_{1}$ successively $m$-times. So $%
J_{(0)}T=T,$ $J_{(1)}T=J_{1}T$ but $J_{m}T\rightarrow M$ is a subbundle of $%
J_{(m)}T\rightarrow M$ for $m\geq 2.$ If $J_{(m)}T\rightarrow M$ admits a
flat linear connection for some $m$, then so does $J_{(k)}T\rightarrow M$
for $k\geq m.$ The smallest such integer $\alpha (M),$ if it exists, is
called the Andreotti invariant of $M.$ If $M$ is $m$-flat for some $M,$ then 
$\mathcal{P}^{\ast }(M,T)=0.$ Therefore $\mathcal{P}^{\ast }(M,T)\neq 0$ $%
\Rightarrow \alpha (M)=\infty ,$ like projective spaces. Now [9] shows that $%
\alpha (M)$ is finite for certain lens spaces and makes a detailed study in
this case.

The reason why $m$-flatness forces $\mathcal{P}^{\ast }(M,T)=0$ as stated in
[9] can be shown as follows. The structure group $j_{(m)}G_{1}(n)$ of $%
J_{(m)}T\rightarrow M$ can be reduced to $G_{1}(n)$ because $%
j_{(m)}G_{1}(n)=G_{1}(n)\rtimes E$ for some subgroup $E$ diffeomorphic to
some $\mathbb{R}^{d}.$ This reduction $F\rightarrow M$ is isomorphic to the
direct sums of certain tensor products of $T$ and $T^{\ast }.$ For $m=2,$
for instance, $F=$ $T\oplus \left( T\otimes T^{\ast }\right) \oplus \left(
T\otimes T^{\ast }\otimes T^{\ast }\right) $ which is easily checked by the
chain rule. Since $J_{(m)}T\rightarrow M$ and $F\rightarrow M$ are
isomorphic we have $\mathcal{P}^{\ast }(M,J_{(m)}T)=\mathcal{P}^{\ast
}(M,F). $ Observe that even though $J_{(m)}T\rightarrow M$ is a natural
bundle of order $m+1,$ $F\rightarrow M$ is a tensor bundle and therefore $%
\mathcal{P}^{\ast }$ is sensitive only to first order jets, a fact which
will be of great importance in Section 7. Therefore, if $M$ is $m$-flat,
then $\mathcal{P}^{\ast }(M,F)=0.$ However the $P$-classes of direct sums
and tensor products are determined by the $P$-classes of the factors
seperately. It follows that the $P$-classes of $F\rightarrow M$, which all
vanish, can be expressed in terms of the $P$-classes of $T\rightarrow M$
which gives polynomial relations of the form

\begin{eqnarray}
ap_{1}(T) &=&0 \\
bp_{1}(T)^{2}+cp_{2}(T) &=&0  \notag \\
dp_{1}(T)^{3}+ep_{2}(T)^{2}+fp_{3}(T) &=&0  \notag \\
...... &=&0  \notag
\end{eqnarray}%
for some constants $a,b,....$ It remains to show $a\neq 0,$ $c\neq 0,$ $%
f\neq 0$ and this is done by explicit algebraic computation. This argument
shows that $m$-flatness for some $m$ implies $\mathcal{P}^{\ast
}(M,J_{(k)}T)=0$ for all $k\geq 0.$

Now [9] shows that $\alpha (M)$ can be arbitrarily large. The following
question will be our driving force in Section 7.

\textbf{Q: }Suppose $\alpha (M)$ is finite so that $\mathcal{P}^{\ast
}(M,J_{(k)}T)=0$ for all $k\geq 0.$ What are the obstructions to $\left(
\alpha (M)-1\right) $-flatness?

Now consider some $\varepsilon \mathcal{G}_{k}$ and the principal $%
\varepsilon \mathcal{G}_{k}^{e,\bullet }\rightarrow M$ with $\mathcal{G}%
_{k}^{e,e}$ as the structure group which we assume to be connected. By the
Iwasawa-Malcev theorem, we decompose $\mathcal{G}_{k}^{e,e}=KE$ where $%
K\subset \mathcal{G}_{k}^{e,e}$ is a maximal compact subgroup and the subset 
$E$ is euclidean.

\begin{lemma}
The restriction of the projection $\pi _{k,1}:\mathcal{G}_{k}^{e,e}%
\rightarrow \mathcal{G}_{1}^{e,e}$ to $K$ is an imbedding.
\end{lemma}

The reason is that the kernel of $\pi _{k,k-1}:\mathcal{G}%
_{k}^{e,e}\rightarrow \mathcal{G}_{k-1}^{e,e}$ is a subgroup of the vector
group $K_{k,k-1}$ in (15) and $K$ must intersect this kernel trivially since
it is compact. Iterating this argument, we see that $K$ must be contained in 
$\mathcal{G}_{1}^{e,e}.$

Let $\mathfrak{\varepsilon G}_{k}\rightarrow M$ be the algebroid of $%
\varepsilon \mathcal{G}_{k}$ and let $P\rightarrow M$ be the principal
bundle associated with $\mathfrak{\varepsilon G}_{k}\rightarrow M$ with $%
GL(m,\mathbb{R})$ as the structure group where $m$ is the dimension of the
fibers of $\mathfrak{\varepsilon G}_{k}\rightarrow M.$ Now $\varepsilon 
\mathcal{G}_{k}^{e,\bullet }\rightarrow M$ is a reduction of $P\rightarrow M$
with $\mathcal{G}_{1}^{e,e}\subset GL(m,\mathbb{R}).$ Using Lemma 22, we can
reduce the structure group further to $K\subset \mathcal{G}_{1}^{e,e}$. Now
we have 
\begin{equation}
\begin{array}{ccc}
I_{GL(m,\mathbb{R})} & \longrightarrow & H_{dR}^{\ast }(M,\mathbb{R)} \\ 
\downarrow \theta &  & \parallel \\ 
I_{K} & \longrightarrow & H_{dR}^{\ast }(M,\mathbb{R}%
\end{array}%
\end{equation}%
where $\theta $ is the restriction homomorphism induced by the restrictions $%
K\subset \mathcal{G}_{k}^{e,e}\subset GL(m,\mathbb{R}).$ Now if $\varepsilon 
\mathcal{G}_{k}$ is locally solvable, then the image of the bottom
homomorphism in (49) vanishes and as above, we believe that this brings
polynomial relations like (48) which depend on $\varepsilon \mathcal{G}_{k}.$
We believe that the clarification of this scenario will explain many known
vanishing phenomena, like Bott vanishing theorem for plane fields, Chern
vanishing theorem which states $\mathcal{P}^{\ast }(M,T)=0$ for a Riemannian
structure of constant curvature, Borel-Hirzebruch vanishing theorem which
states $\mathcal{P}^{\ast }(M,T)=0$ for $G/T$ where $G$ is compact and $T$
is a maximal torus, the well known relations between the Chern classes of
projective spaces...and many other phenomena about the structure of the
characteristic classes of homogeneous spaces.

\section{Higher order characteristic classes}

This section is the main core of this note. We will outline here the
construction of higher order obstructions to local solvability. In
particular, our method will give obstructions to $m$-flatness. We do not
know whether these invariants can be nontrivial. What we do know, however,
is that they will be highly nontrivial if they are nontrivial at all!

As we observed in Section 6, the $P$-algebra algebra $\mathcal{P}^{\ast }(M,%
\mathfrak{G}_{k})\subset H_{dR}^{\ast }(M,\mathbb{R})$ is sensitive only to
first order jets and is "topological" even though $k$ is large. This
topology persists even if $\varepsilon \mathcal{G}_{k}$ is locally solvable
in the following sense.

Let a connected $G$ act transitively on $M.$ If $M$ is compact and the
stabilizers of the action are connected (this is so if $M$ is also simply
connected), then a maximal compact subgroup $K\subset G$ acts also
transitively on $M$ according to [13]. Hence $M=G/H=K/K\cap H$ for some
stabilizer $H.$ Therefore, as long as we are interested in the topological
properties of a compact Klein geometry $G/H$ with connected $H$, we may
assume that $G$ (and therefore $H$) is compact. Obviously $ord(G/H)=0$ if $%
H=\{1\}.$ Now Lemma 22 implies

\begin{proposition}
If $H$ is compact, connected and nontrivial, then $ord(G/H)=1.$
\end{proposition}

The above arguments make it clear that the cohomological invariants of $%
\varepsilon \mathcal{G}_{k}$ which depend on $k$ can not be topological if $%
k $ is large. Further, we should not search for such invariants in the
cohomology of the base $M$. In particular, we should not consider $%
\varepsilon \mathcal{G}_{k}$ as a fibering over $M$ as we did so far, i.e.,
we should not let $\varepsilon \mathcal{G}_{k}$ act on $M.$

Inspired by Proposition 14, we start with the following

\begin{proposition}
Let $\varepsilon \mathcal{G}_{k}$ be a $phg$ of order $k$ on $M.$ Then the
the total space $\mathcal{G}_{k}^{e,\bullet }$ of the principal bundle $%
\mathcal{G}_{k}^{e,\bullet }\rightarrow M$ is parallelized by $\varepsilon $
in a canonical way.
\end{proposition}

Proposition 24 is trivial for $k=0$ since $\mathcal{G}_{0}^{e,\bullet
}\simeq M$ is parallelizable by definition. If $\mathcal{R}_{k}\mathcal{=}0$
and $\widetilde{\varepsilon \mathcal{G}_{k}}$ is globalizable, then $%
\mathcal{G}_{k}^{e,\bullet }\simeq G$ which is surely parallelizable. The
surprising fact is that the statement holds for all $k\geq 0$ without the
assumption $\mathcal{R}_{k}=0.$

The proof follows almost trivially from the definition of a $phg$. First we
recall the following trivial

Fact: Let $M$ be a smooth manifold and $p,q\in M.$ There is a canonical
identification between the following sets:

$i)$ The set of $1$-arrows from $p$ to $q$

$ii)$ The set of isomorphisms $T_{p}(M)\rightarrow T_{q}(M)$

Now let $\overline{p},\overline{q}$ be arbitrary points on $\mathcal{G}%
_{k}^{e,\bullet }$ which project to $p,q\in M$. So $\overline{p},\overline{q}
$ are two $k$-arrows from $e$ to $p,q$ respectively, say $\overline{p}%
=j_{k}(f)^{e,p}$ and $\overline{q}=j_{k}(g)^{e,q}.$ Therefore $%
j_{k}(g)^{e,q}\circ \left[ j_{k}(f)^{e,p}\right] ^{-1}=$ $j_{k}\left( g\circ
f^{-1}\right) ^{p,q}$ is a $k$-arrow from $p$ to $q.$ The splitting $%
\varepsilon $ gives the \textit{unique }$\left( k+1\right) $-arrow $%
\varepsilon j_{k}\left( g\circ f^{-1}\right) ^{p,q}$ from $p$ to $q.$ Now
let $\xi _{\overline{p}}$ be a tangent vector at $\overline{p}.$ According
to (22) the tangent space $T(\mathcal{G}_{k}^{e,\bullet })^{\overline{p}}$
is the same as the fiber $\mathfrak{G}_{k}^{p}$ of the algebroid $\mathfrak{G%
}_{k}\rightarrow M$ over $p.$ Therefore the isomorphism $\varepsilon
j_{k}\left( g\circ f^{-1}\right) ^{p,q}:\mathfrak{G}_{k}^{p}\rightarrow 
\mathfrak{G}_{k}^{q}$ in (26) is an isomorphism $\varepsilon j_{k}\left(
g\circ f^{-1}\right) ^{p,q}:$ $T(\mathcal{G}_{k}^{e,\bullet })^{\overline{p}%
}\rightarrow T(\mathcal{G}_{k}^{e,\bullet })^{\overline{q}}$. The above fact
implies that the object $\varepsilon j_{k}\left( g\circ f^{-1}\right) ^{p,q}$%
, which is a $(k+1)$- arrow from $p$ to $q$ is at the same time a $1$-arrow
from $\overline{p}$ to $\overline{q}.$ In short, $\left( k+1\right) $-arrows
on $M$ define $1$-arrows on $\mathcal{G}_{k}^{e,\bullet }!!$ So for any two
points $\overline{p},\overline{q}\in \mathcal{G}_{k}^{e,\bullet }$ we
established a unique $1$-arrow from $\overline{p}$ to $\overline{q}$ which
is equivalent to the parallelizability of $\mathcal{G}_{k}^{e,\bullet }$
since this assignment is smooth and is a homomorphism of groupoids. We will
continue to denote this splitting by $\varepsilon $.

The above proof warns us that we should be more careful with our notation.
So we denote \ a $phg$ $\varepsilon \mathcal{G}_{k}$ on $M$ by $(\mathcal{G}%
_{0}^{e,\bullet },\varepsilon \mathcal{G}_{k})=(M,\varepsilon \mathcal{G}%
_{k})$ henceforth. Accordingly we will denote the $phg$ of order zero in
Proposition 24 by $(\mathcal{G}_{k}^{e,\bullet },\varepsilon \mathcal{G}%
_{k}).$ Let $\mathcal{R}_{0}$ denote the curvature of $(\mathcal{G}%
_{k}^{e,\bullet },\varepsilon \mathcal{G}_{k}).$ By Propositions 17,19, $%
\mathcal{R}_{0}=0\Leftrightarrow $ $(\mathcal{G}_{k}^{e,\bullet
},\varepsilon \mathcal{G}_{k})$ integrates to the pseudogroup $\widetilde{(%
\mathcal{G}_{k}^{e,\bullet },\varepsilon \mathcal{G}_{k})}$, i.e., all $1$%
-arrows of $(\mathcal{G}_{k}^{e,\bullet },\varepsilon \mathcal{G}_{k})$
integrate uniquely to local diffeomorphisms on $\mathcal{G}_{k}^{e,\bullet
}. $ If $\widetilde{(\mathcal{G}_{k}^{e,\bullet },\varepsilon \mathcal{G}%
_{k})}$ is globalizable, then we get the transformation group $G$ which acts 
\textit{simply transitively }on $\mathcal{G}_{k}^{e,\bullet }$ and $\mathcal{%
G}_{k}^{e,\bullet }\simeq G$ as before. It is extremely crucial to observe
that the action of $G$ may not descend to $M$, i.e., $G$ may not act on $M:$
for this we need the stronger condition $\mathcal{R}_{k}=0$ which implies $%
\mathcal{R}_{0}=0.$

Now the algebroid of $(\mathcal{G}_{k}^{e,\bullet },\varepsilon \mathcal{G}%
_{k})$ is $\varepsilon \mathfrak{G}_{k}\rightarrow \mathcal{G}%
_{k}^{e,\bullet }$ which is simply $T(\mathcal{G}_{k}^{e,\bullet
})\rightarrow \mathcal{G}_{k}^{e,\bullet }$ together with the splitting $%
\varepsilon :T(\mathcal{G}_{k}^{e,\bullet })\rightarrow \varepsilon T(%
\mathcal{G}_{k}^{e,\bullet })\subset J_{1}T(\mathcal{G}_{k}^{e,\bullet }).$
We will denote this algebroid by $(\mathcal{G}_{k}^{e,\bullet },\varepsilon 
\mathfrak{G}_{k})$ for notational convinience below. Let $\mathfrak{R}_{0}$
be the curvature of $(\mathcal{G}_{k}^{e,\bullet },\varepsilon \mathfrak{G}%
_{k})$ obtained by linearizing $\mathcal{R}_{0}.$ For $\overline{p}\in 
\mathcal{G}_{k}^{e,\bullet }$ and $\xi _{p},\eta _{p}\in \mathfrak{G}%
_{k}^{p}=T(\mathcal{G}_{k}^{e,\bullet })^{\overline{p}}$, we have

\begin{equation}
\mathfrak{R}_{0}(\overline{p})(\xi _{p},\eta _{p})\in Hom(\mathfrak{G}%
_{k}^{p},\mathfrak{G}_{k}^{p})
\end{equation}

Clearly the $P$-algebra

\begin{equation}
\mathcal{P}^{\ast }(\mathcal{G}_{k}^{e,\bullet },\varepsilon \mathfrak{G}%
_{k})\subset H_{dR}^{\ast }(\mathcal{G}_{k}^{e,\bullet },\mathbb{R})
\end{equation}%
vanishes since $\mathcal{G}_{k}^{e,\bullet }$ is parallelizable. We define
the $2i$-forms $Tr(\mathfrak{R}_{0}^{i})$ on $\mathcal{G}_{k}^{e,\bullet }$
by

\begin{equation}
Tr(\mathfrak{R}_{0}^{i})(\overline{p},\xi _{p}^{1},\xi _{p}^{2},...,\xi
_{p}^{2i})\overset{def}{=}\frac{1}{(2k)!}\dsum\limits_{\sigma }sgn(\sigma
)\left( Tr\left( \mathfrak{R}_{0}(\overline{p})(\xi _{p}^{1},\xi
_{p}^{2})\circ ....\circ \mathfrak{R}_{0}(\overline{p})(\xi _{p}^{2i-1},\xi
_{p}^{2i}\right) \right)
\end{equation}%
where the summation is taken over all permutations $\sigma $ of $%
(1,2,...,2j).$ The forms $Tr(\mathfrak{R}_{0}^{i})$ are exact in the de Rham
complex of $\mathcal{G}_{k}^{e,\bullet }.$ In fact, the "Chern-Simons forms"
with a surprisingly different interpretation supply some canonical
primitives for $Tr(\mathfrak{R}_{k}^{i})$ (see Section 10).

Now we have the following crucial

\begin{lemma}
$Tr(\mathfrak{R}_{0}^{i})(\overline{p}\circ a)=Tr(\mathfrak{R}_{0}^{i})(%
\overline{p})$ for all $a\in \mathcal{G}_{k}^{e,e}$ and $\overline{p}\in 
\mathcal{G}_{k}^{e,\bullet }.$
\end{lemma}

Lemma 25 states that $Tr(\mathfrak{R}_{0}^{i})$ is a right invariant $2i$%
-form on $\mathcal{G}_{k}^{e,\bullet }\rightarrow M.$ Observe that $Tr(%
\mathfrak{R}_{k}^{i})$ in Section 6 which generate $\mathcal{P}^{\ast }(M,%
\mathfrak{G}_{k})\subset H_{dR}^{\ast }(M,\mathbb{R})$ also live on $%
\mathcal{G}_{k}^{e,\bullet }$ but they are horizontal over $M$ and therefore
descend to $M.$ However the forms $Tr(\mathfrak{R}_{0}^{i})$ are surely not
horizontal over $M$ unless $k=0.$

Now let $\wedge ^{k}(\mathcal{G}_{k}^{e,\bullet },M)$ denote the space of $k$%
-forms on $\mathcal{G}_{k}^{e,\bullet }$ which are right invariant over $M.$
The exterior derivative $d$ of the de Rham complex of $\mathcal{G}%
_{k}^{e,\bullet }$ restricts as $d:$ $\wedge ^{k}(\mathcal{G}_{k}^{e,\bullet
},M)\rightarrow \wedge ^{k+1}(\mathcal{G}_{k}^{e,\bullet },M)$ and we have
the subcomplex

\begin{equation}
\wedge ^{0}(\mathcal{G}_{k}^{e,\bullet },M)\overset{d}{\longrightarrow }%
\wedge ^{1}(\mathcal{G}_{k}^{e,\bullet },M)\overset{d}{\longrightarrow }...%
\overset{d}{\longrightarrow }\wedge ^{s}(\mathcal{G}_{k}^{e,\bullet },M)
\end{equation}%
where $s=\dim \mathcal{G}_{k}^{e,\bullet }$ and $\wedge ^{0}(\mathcal{G}%
_{k}^{e,\bullet },M)=C^{\infty }(\mathcal{G}_{k}^{e,\bullet }).$ The
cohomology $H_{inv}^{\ast }(\mathcal{G}_{k}^{e,\bullet },M)$ of (53) is
called the algebroid cohomology of $\mathfrak{G}_{k}\rightarrow M$ which we
write also as $H^{\ast }(M,\mathfrak{G}_{k}).$ For $k=0$ (53) is the de Rham
complex of $M$ but for $k\geq 1$ it is a \textit{proper }subcomplex.

The crucial fact now is that the forms $Tr\left( \mathfrak{R}_{0}^{i}\right) 
$ which are exact in the de Rham complex of $\mathcal{G}_{k}^{e,\bullet }$
need not be exact in (53) for $k\geq 1.$ Therefore the forms $Tr\left( 
\mathfrak{R}_{0}^{i}\right) $ generate a subalgebra $\widehat{\mathcal{P}%
^{\ast }}(M,\varepsilon \mathfrak{G}_{k})\subset H_{inv}^{2i}(\mathcal{G}%
_{k}^{e,\bullet },M)=$ $H^{\ast }(M,\mathfrak{G}_{k})$

\begin{definition}
The algebra $\widehat{\mathcal{P}^{\ast }}(M,\varepsilon \mathfrak{G}%
_{k})\subset H_{inv}^{2i}(\mathcal{G}_{k}^{e,\bullet },M)=$ $H^{\ast }(M,%
\mathfrak{G}_{k})$ is the $k$-th order Pontryagin algebra of the $phg$ $%
\varepsilon \mathcal{G}_{k}.$\textit{\ }
\end{definition}

Observe that $\widehat{\mathcal{P}^{\ast }}(M,\varepsilon \mathfrak{G}_{k})$
depends on $\varepsilon $ which is fixed by the definition of $\varepsilon 
\mathcal{G}_{k}.$ Clearly $\widehat{\mathcal{P}^{\ast }}(M,\varepsilon 
\mathfrak{G}_{k})=0$ for $k=0$ since $M$ is parallelizable.

The next proposition shows that the above construction of $(\mathcal{G}%
_{k}^{e,\bullet },\varepsilon \mathcal{G}_{k})$ is a particular case.

\begin{proposition}
The pair $(\mathcal{G}_{i}^{e,\bullet },\varepsilon \mathcal{G}_{k})$
defines a $phg$ of order $k-i$ on the total space $\mathcal{G}%
_{i}^{e,\bullet }$ of the principal bundle $\mathcal{G}_{i}^{e,\bullet
}\rightarrow M$ for $0\leq i\leq k.$
\end{proposition}

The main idea of Proposition 27 is simple: Let $f$ be a local diffeomorphism
on $M$ with $f(p)=q$ and $j_{i}(f)^{p,q}\in \mathcal{G}_{i}^{p,q}$. Now $f$
defines a function $\mathcal{G}_{i}^{e,p}\rightarrow \mathcal{G}_{i}^{e,q}$
by $\alpha _{i}^{e,p}\rightarrow $ $j_{i}(f)^{p,q}\circ \alpha _{i}^{e,p}.$
If $j_{i+1}(f)^{p,q}\in \mathcal{G}_{i+1}^{p,q}$ then $j_{i+1}(f)^{p,q}$
defines a $1$-arrow from $\alpha _{i}^{e,p}$ to $j_{i}(f)^{p,q}\circ \alpha
_{i}^{e,p}$ for any $\alpha _{i}^{e,p}\in \mathcal{G}_{i}^{e,p}.$ Similarly $%
j_{i+2}(f)^{p,q}\in \mathcal{G}_{i+2}^{p,q}$ defines a $2$-arrow from $%
\alpha _{i}^{e,p}$ to $j_{i}(f)^{p,q}\circ \alpha _{i}^{e,p}.$ Iterating
this process we see that $j_{k}(f)^{p,q}\in \mathcal{G}_{k}^{p,q}$ defines a 
$(k-i)$-arrow from $\alpha _{i}^{e,p}$ to $j_{i}(f)^{p,q}\circ \alpha
_{i}^{e,p}.$ Finally, above any such $(k-i)$-arrow, there is a unique $%
(k-i+1)$-arrow defined by $\varepsilon j_{k}(f)^{p,q}\in \varepsilon 
\mathcal{G}_{k}^{p,q}.$

Let $\mathcal{R}_{k-i}$ denote the curvature of $(\mathcal{G}_{i}^{e,\bullet
},\varepsilon \mathcal{G}_{k})$. Now $\mathcal{R}_{k-i}=0\Leftrightarrow $
the pseudogroup $\widetilde{(\varepsilon \mathcal{G}_{k},\mathcal{G}%
_{k}^{e,\bullet })}$ which acts on $\mathcal{G}_{k}^{e,\bullet }$ descends
to the pseudogroup $\widetilde{(\mathcal{G}_{i}^{e,\bullet },\varepsilon 
\mathcal{G}_{k},)}$ which acts on $\mathcal{G}_{i}^{e,\bullet }.$ In
particular, we have

\begin{proposition}
$\mathcal{R}_{k-i}=0\Rightarrow \mathcal{R}_{k-i-1}=...=\mathcal{R}_{0}=0$
\end{proposition}

The algebroid $(\mathcal{G}_{i}^{e,\bullet },\varepsilon \mathfrak{G}_{k})$
of $(\mathcal{G}_{i}^{e,\bullet },\varepsilon \mathcal{G}_{k})$ is easily
described explicitly: $\mathfrak{G}_{i}^{p}$ is the tangent space at $%
\overline{p}\in \mathcal{G}_{i}^{e,p}$, $\mathfrak{G}_{i+1}^{p}$ is $1$-jets
of vector fields at $\overline{p}$....and $\mathfrak{G}_{k}^{p}$ is the $%
(k-i)$-jets of vector fields at $\overline{p}.$ Finally, above any such $%
(k-i)$-jet, there is a unique $(k-i+1)$-jet defined by $\varepsilon 
\mathfrak{G}_{k}^{p}.$

Let $\mathfrak{R}_{k-i}$ be the curvature of $(\mathcal{G}_{i}^{e,\bullet
},\varepsilon \mathfrak{G}_{k}).$ We have

\begin{equation}
\mathfrak{R}_{k-i}(\overline{p},\xi _{p},\eta _{p})\in Hom(\mathfrak{G}%
_{k}^{p},\mathfrak{G}_{k}^{p})\text{ \ \ \ }\overline{p}\in \mathcal{G}%
_{i}^{e,\bullet },\text{ }\xi _{p},\eta _{p}\in \mathfrak{G}_{i}^{p}
\end{equation}%
where $\overline{p}\in \mathcal{G}_{i}^{e,\bullet }$ projects to $p\in M$
and $\xi _{p},\eta _{p}\in \mathfrak{G}_{i}^{p}=T(\mathcal{G}_{i}^{e,\bullet
})^{\overline{p}}.$ The algebra $\widehat{\mathcal{P}^{\ast }}(\mathcal{G}%
_{0}^{e,\bullet },\varepsilon \mathfrak{G}_{k})=$ $\mathcal{P}^{\ast }(%
\mathcal{G}_{0}^{e,\bullet }\varepsilon \mathfrak{G}_{k})$ is considered in
Section 6 and we defined $\widehat{\mathcal{P}^{\ast }}(\mathcal{G}%
_{k}^{e,\bullet },\varepsilon \mathfrak{G}_{k})$ above. Henceforth we assume 
$1\leq i\leq k-1.$

Now using (54) and (52) we define the forms $Tr\mathfrak{R}_{k-i}^{j}$ in
the de Rham complex of $\mathcal{G}_{i}^{e,\bullet }$ which are all exact
because the principal bundle

\begin{equation}
\mathcal{G}_{k}^{e,\bullet }\rightarrow \mathcal{G}_{i}^{e,\bullet }
\end{equation}%
has contractible fibers and is therefore trivial. However, as in Lemma 25,
we have

\begin{equation}
\mathfrak{R}_{k-i}(\overline{p}a)=\mathfrak{R}_{k-i}(\overline{p})
\end{equation}%
where $a\in \mathcal{G}_{i}^{e,e}$ and we consider

\begin{equation}
\left[ Tr\mathfrak{R}_{k-i}^{j}\right] \in H_{inv}^{2j}(\mathcal{G}%
_{i}^{e,\bullet },M)=H_{inv}^{\ast }(\mathcal{G}_{i}^{e,\bullet },M)
\end{equation}

\begin{definition}
The subalgebra $\widehat{\mathcal{P}^{\ast }}(\mathcal{G}_{i}^{e,\bullet
},\varepsilon \mathfrak{G}_{k})\subset H_{inv}^{\ast }(\mathcal{G}%
_{i}^{e,\bullet },M)$ generated by the forms (54) is the $i$-th order
Pontryagin algebra of $\varepsilon \mathcal{G}_{k},$ $1\leq i\leq $ $k-1.$
\end{definition}

Obviously $\widehat{\mathcal{P}^{\ast }}(\mathcal{G}_{i}^{e,\bullet
},\varepsilon \mathfrak{G}_{k})=0$ if $(\mathcal{G}_{i}^{e,\bullet
},\varepsilon \mathfrak{G}_{k})$ is locally solvable.

Observe that there seems to be no reason for $\left[ Tr\mathfrak{R}_{k-i}^{j}%
\right] =0$ for $j$ odd unless $i=0.$

Unfortunately we do not get new invariants in Riemannian geometry because $%
H_{inv}^{\ast }(\mathcal{G}_{1}^{e,\bullet },M)=H_{dR}^{\ast }(\mathcal{G}%
_{1}^{e,\bullet },\mathbb{R})$ since $O(n)$ is compact.

The above method gives also obstructions to $m$-flatness as follows.

We have the groupoids $\mathcal{U}_{(k)}\overset{def}{=}J_{(k)}(M\times M)$
with the algebroids $J_{(k)}T\rightarrow M,$ $k\geq 0.$ There is a 1-1
correspondence between the following objects.

$i)$ Sections of $J_{(k+1)}T\rightarrow M$

$ii)$ Connections on the vector bundle $J_{(k)}T\rightarrow M$

Any such object $\varepsilon _{k}$ defines a connection (using the same
notation) $\varepsilon _{k}$ on the tangent bundle $T\mathcal{U}%
_{(k)}^{e,\bullet }\rightarrow \mathcal{U}_{(k)}^{e,\bullet }$ where $%
\mathcal{U}_{(k)}^{e,\bullet }\rightarrow M$ is the principal bundle of the
groupoid $\mathcal{U}_{(k)}$ with base point $e\in M.$ We apply the
Chern-Weil construction to $\varepsilon _{k}$ and get the subalgebra $%
\mathcal{P}^{\ast }(\mathcal{U}_{(k)}^{e,\bullet },T\mathcal{U}%
_{(k)}^{e,\bullet })\subset H_{dR}^{\ast }(\mathcal{U}_{(k)}^{e,\bullet },%
\mathbb{R)}$ which is trivial. The forms obtained in this way are right
invariant on $\mathcal{U}_{(k)}^{e,\bullet }\rightarrow M$ and a change of $%
\varepsilon _{k}$ adds to such a form a boundary which is right invariant
(compare to Proposition 31 below). Thus we get the subalgebra $\widehat{%
\mathcal{P}^{\ast }}(\mathcal{U}_{(k)}^{e,\bullet },T\mathcal{U}%
_{(k)}^{e,\bullet })$ $\subset H_{inv}^{\ast }(\mathcal{U}_{(k)}^{e,\bullet
},M)$ $=$ $H^{\ast }(M,J_{(k)}T)$ = the cohomology of the algebroid $%
J_{(k)}T\rightarrow M.$ Clearly $\widehat{\mathcal{P}^{\ast }}(\mathcal{U}%
_{(k)}^{e,\bullet },T\mathcal{U}_{(k)}^{e,\bullet })=0$ if $M$ is $k$-flat.
This construction is the analog of Proposition 24.

More generally, we define 
\begin{equation}
\widehat{\mathcal{P}^{\ast }}(\mathcal{U}_{(r)}^{e,\bullet },J_{(k-r)}T%
\mathcal{U}_{(r)}^{e,\bullet })\subset H_{inv}^{\ast }(\mathcal{U}%
_{(r)}^{e,\bullet },M)=H^{\ast }(J_{(r)}T,M)
\end{equation}
for $0\leq r\leq k$ by observing that $\varepsilon _{k}$ gives a connection
on the vector bundle $J_{(k-r)}T\mathcal{U}_{(r)}^{e,\bullet }$ $\rightarrow 
\mathcal{U}_{(r)}^{e,\bullet },$ $0\leq r\leq k.$ For $r=0$ we get the
topological obstructions $\widehat{\mathcal{P}^{\ast }}(M,J_{(k)}T)=\mathcal{%
P}^{\ast }(M,J_{(k)}T)\subset H_{dR}^{\ast }(M,\mathbb{R}).$

To summarize, we have

\begin{proposition}
The restricted $P$-algebras $\widehat{\mathcal{P}^{\ast }}(\mathcal{U}%
_{(r)}^{e,\bullet },J_{(k-r)}T\mathcal{U}_{(r)}^{e,\bullet })\subset
H_{inv}^{\ast }(\mathcal{U}_{(r)}^{e,\bullet },M)$ $=H^{\ast }(J_{(r)}T,M)$
for $0\leq r\leq k$ vanish if $M$ is $k$-flat.
\end{proposition}

\section{ Dependence of the isomorphism class}

Up to now we dealt with some fixed $\varepsilon \mathcal{G}_{k}.$ Now we
want to define the notion of equivalence of $phg$'s in such a way that the
above constructions depend only on the equivalence class of $\varepsilon 
\mathcal{G}_{k}.$ As a crucial point, we will not fix the principal bunde $%
\mathcal{G}_{k}^{e,\bullet }\rightarrow M$ and change the connection $%
\varepsilon $ but change $\varepsilon \mathcal{G}_{k}$ and preserve the
order of jets and the vertex class $\left\{ \varepsilon \mathcal{G}%
_{k}\right\} $.

So we start with some $\varepsilon \mathcal{G}_{k}$ on $M.$ Consider the
group bundle $\mathcal{A}_{k+1}\overset{def}{=}\cup _{x\in M}\mathcal{U}%
_{k+1}^{x,x}$ $\rightarrow M.$ A smooth section of $\mathcal{A}%
_{k+1}\rightarrow M$ is called a gauge transformation of order $k+1$ on $M.$
The set $\Gamma \mathcal{A}_{k+1}$ of gauge transformations is a group with
fiberwise composition. Now $a\in \Gamma \mathcal{A}_{k+1}$ acts on the
arrows of $\varepsilon \mathcal{G}_{k}$ by

\begin{equation}
a\cdot \varepsilon j_{k}(f)^{p,q}\overset{def}{=}a(q)\circ \varepsilon
j_{k}(f)^{p,q}\circ a(p)^{-1}
\end{equation}

We see that $a\cdot \varepsilon \mathcal{G}_{k}$ is another $phg$ having the
same vertex class as $\varepsilon \mathcal{G}_{k}$, i.e., $\left\{ a\cdot
\varepsilon \mathcal{G}_{k}\right\} =\left\{ \varepsilon \mathcal{G}%
_{k}\right\} .$ We have the projections $\pi :\mathcal{A}_{k+1}\rightarrow 
\mathcal{A}_{i},$ $1\leq i\leq k,$ which give the projections $\pi :\Gamma 
\mathcal{A}_{k+1}\rightarrow \Gamma \mathcal{A}_{i}$. By projecting (59) to
the jets of order $i,$ we obtain the commutative diagram

\begin{equation}
\begin{array}{ccc}
\mathcal{G}_{k}^{e,\bullet } & \rightarrow & a\cdot \mathcal{G}%
_{k}^{e,\bullet } \\ 
\downarrow \pi &  & \downarrow \pi \\ 
\mathcal{G}_{i}^{e,\bullet } & \rightarrow & \pi (a)\cdot \mathcal{G}%
_{i}^{e,\bullet }%
\end{array}%
\end{equation}%
where the horizontal arrows are isomorphisms of principal bundles.

The action of $\Gamma \mathcal{A}_{k+1}$ on $\varepsilon \mathcal{G}_{k}$
gives a natural action of $\Gamma \mathcal{A}_{k+1}$ on the algebroid $%
\varepsilon \mathfrak{G}_{k}.$ We have the commutative diagram

\begin{equation}
\begin{array}{ccc}
\varepsilon \mathcal{G}_{k} & \Longrightarrow & \varepsilon \mathfrak{G}_{k}
\\ 
\downarrow &  & \downarrow \\ 
a\cdot \varepsilon \mathcal{G}_{k} & \Longrightarrow & a\cdot \varepsilon 
\mathfrak{G}_{k}%
\end{array}%
\end{equation}%
where $\Longrightarrow $ denotes linearization. We write $\varepsilon 
\mathcal{G}_{k}\sim $ $\varepsilon ^{\prime }\mathcal{G}_{k}^{\prime }$ if $%
\varepsilon ^{\prime }\mathcal{G}_{k}^{\prime }=a\cdot \varepsilon \mathcal{G%
}_{k}$ for some $a\in \Gamma \mathcal{A}_{k+1}.$ We have $\varepsilon 
\mathcal{G}_{k}\sim $ $\varepsilon ^{\prime }\mathcal{G}_{k}^{\prime
}\Leftrightarrow \varepsilon \mathfrak{G}_{k}\sim $ $\varepsilon ^{\prime }%
\mathfrak{G}_{k}^{\prime }$. We denote the equivalence classes of $%
\varepsilon \mathcal{G}_{k},$ $\varepsilon \mathfrak{G}_{k}$ by $\left[
\varepsilon \mathcal{G}_{k}\right] ,$ $\left[ \varepsilon \mathfrak{G}_{k}%
\right] .$ Clearly the vector bundles $\varepsilon \mathfrak{G}%
_{k}\rightarrow M$ and $a\cdot \varepsilon \mathfrak{G}_{k}\rightarrow M$
are isomorphic. However the $phg$'s $\varepsilon \mathfrak{G}_{k}\rightarrow
M$ and $\varepsilon ^{\prime }\mathfrak{G}_{k}^{\prime }\rightarrow M$ may
be isomorphic as vector bundles but inequivalent as $phg$'s as defined
above. The main point is that the above equivalence respects the order of
jets whereas the topological concept of "vector bundle isomorphism" does not.

Now the assignment $\varepsilon \mathfrak{G}_{k}\Rightarrow \mathcal{P}%
^{\ast }(\mathcal{G}_{0}^{e,\bullet },\varepsilon \mathfrak{G}_{k})$ is
rather crude as it depends on the isomorphism class of the vector bundle $%
\varepsilon \mathfrak{G}_{k}\rightarrow M.$ However, it turns out that the
assignments $\varepsilon \mathfrak{G}_{k}\Rightarrow $ $\widehat{\mathcal{P}%
^{\ast }}(\mathcal{G}_{i}^{e,\bullet },\varepsilon \mathfrak{G}_{k}),$ $%
1\leq i\leq k$ depend on $\left[ \varepsilon \mathfrak{G}_{k}\right] $ in a
sense to be made precise below.

Recall that $\mathcal{C}^{\ast }(\mathcal{G}_{i}^{e,\bullet },M)$ denotes
the complex of right invariant forms on the principal bundle $\mathcal{G}%
_{i}^{e,\bullet }\rightarrow M,$ i.e., the complex computing the algebroid
cohomology of $\mathfrak{G}_{i}\rightarrow M.$ The linearization of the
bottom isomorphism of (61) shows that $a\in \Gamma \mathcal{A}_{k+1}$
defines an isomorphism

\begin{equation}
a^{\ast }:\mathcal{C}^{\ast }(\mathcal{G}_{i}^{e,\bullet },M)\rightarrow 
\mathcal{C}^{\ast }(\pi (a)\cdot \mathcal{G}_{i}^{e,\bullet },M)
\end{equation}%
for $0\leq i\leq k$. For $i=0$, $a^{\ast }=Id$ and both complexes are the de
Rham complex of $M.$ It follows that $a^{\ast }$ acts as an isomorphism on
the complex (53) as

\begin{equation}
\begin{array}{ccccccc}
\wedge ^{0}(\mathcal{G}_{i}^{e,\bullet },M) & \overset{d}{\longrightarrow }
& \wedge ^{1}(\mathcal{G}_{i}^{e,\bullet },M) & \overset{d}{\longrightarrow }
& .... & \overset{d}{\longrightarrow } & \wedge ^{s}(\mathcal{G}%
_{i}^{e,\bullet },M) \\ 
\downarrow a^{\ast } &  & \downarrow a^{\ast } &  & \downarrow a^{\ast } & 
& \downarrow a^{\ast } \\ 
\wedge ^{0}(a\cdot \mathcal{G}_{i}^{e,\bullet },M) & \overset{d}{%
\longrightarrow } & \wedge ^{1}(a\cdot \mathcal{G}_{i}^{e,\bullet },M) & 
\overset{d}{\longrightarrow } & .... & \overset{d}{\longrightarrow } & 
\wedge ^{s}(a\cdot \mathcal{G}_{i}^{e,\bullet },M)%
\end{array}%
\end{equation}

Therefore $a^{\ast }$ induces an isomorphism

\begin{equation}
a^{\ast }:H^{\ast }(\mathcal{G}_{i}^{e,\bullet },M)\rightarrow H^{\ast
}(a\cdot \mathcal{G}_{i}^{e,\bullet },M)
\end{equation}

Consider the forms $Tr\mathfrak{R}_{k-i}^{j}\in \wedge ^{2j}(\mathcal{G}%
_{i}^{e,\bullet },M)$. Let $Tr\mathfrak{R}_{k-i}^{j}(a)\in \wedge
^{2j}(a\cdot \mathcal{G}_{i}^{e,\bullet },M)$ denote the forms constructed
on $a\cdot \varepsilon \mathcal{G}_{k}$ using the curvature $\mathfrak{R}%
_{k}(a)$ of $(a\cdot \mathcal{G}_{i}^{e,\bullet },M)$ in the same way. Thus $%
a^{\ast }Tr\mathfrak{R}_{k-i}^{j}$ and $Tr\mathfrak{R}_{k-i}^{j}(a)$ both
live in $\wedge ^{2j}(M,a\cdot \mathcal{G}_{i}^{e,\bullet }).$

Now we have the following fundamental

\begin{proposition}
$a^{\ast }Tr\mathfrak{R}_{k-i}^{j}-Tr\mathfrak{R}_{k-i}^{j}(a)$ is exact in
the bottom complex of (63).
\end{proposition}

\begin{corollary}
$a^{\ast }$ induces an isomorphism
\end{corollary}

\QTP{Body Math}
\begin{equation}
a^{\ast }:\widehat{\mathcal{P}^{\ast }}(\mathcal{G}_{i}^{e,\bullet
},\varepsilon \mathfrak{G}_{k})\longrightarrow \widehat{\mathcal{P}^{\ast }}%
(\pi (a)\cdot \mathcal{G}_{i}^{e,\bullet },a\cdot \varepsilon \mathfrak{G}%
_{k})
\end{equation}

\begin{corollary}
We have the well defined assignments%
\begin{equation}
\left[ \varepsilon \mathcal{G}_{k}\right] \Longrightarrow \widehat{\mathcal{P%
}^{\ast }}(\mathcal{G}_{i}^{e,\bullet },\varepsilon \mathfrak{G}_{k})\text{
\ \ \ \ }0\leq i\leq k
\end{equation}%
In particular, $\widehat{\mathcal{P}^{\ast }}(\mathcal{G}_{i}^{e,\bullet
},\varepsilon \mathfrak{G}_{k})=0,$\ $0\leq i\leq k,$ if the equivalence
class $\left[ \varepsilon \mathcal{G}_{k}\right] $ contains a locally
solvable $phg.$
\end{corollary}

\QTP{Body Math}
If $U\subset M$ is an open subset, we define the restriction $\varepsilon 
\mathcal{G}_{k\mid U}$ of $\varepsilon \mathcal{G}_{k}$ as the arrows whose
source and targets are contained in $U.$ Clearly $\varepsilon \mathcal{G}%
_{k\mid U}$ also satisfies $i),$ $ii)$ of Definition 1 and therefore defines
a $phg$ on $U.$ By Corollary 33 we obtain the algebras $\widehat{\mathcal{P}%
^{\ast }}(\mathcal{G}_{i\mid U}^{e,\bullet },\varepsilon \mathfrak{G}_{k\mid
U})$.

\begin{proposition}
$\widehat{\mathcal{P}^{\ast }}(\mathcal{G}_{i\mid U}^{e,\bullet
},\varepsilon \mathfrak{G}_{k\mid U})=0$ for a coordinate neighborhood $%
U\subset M$ diffeomorphic to $\mathbb{R}^{n}.$
\end{proposition}

\QTP{Body Math}
To see this, $iii)$ of Definition 1 implies $\left\langle \varepsilon 
\mathcal{G}_{k}\right\rangle =\left\langle H\right\rangle _{G}$ for some
homogeneous space $G/H=N.$ Let $\varepsilon \mathcal{G}_{k}^{\prime }$ be
the globally solvable $phg$ on $N$ defined by $G/H$ and consider $%
\varepsilon ^{\prime }\mathcal{G}_{k\mid V}^{\prime }$ for some $V\subset N$
diffeomorphic to $\mathbb{R}^{n}$ which gives an identification $U\simeq V.$
So we have the two $phg$'s $\varepsilon \mathcal{G}_{k\mid V}^{\prime }$ and 
$\varepsilon ^{\prime }\mathcal{G}_{k\mid V}$ defined on $V$ and $%
\varepsilon ^{\prime }\mathcal{G}_{k\mid V}^{\prime }$ is locally solvable.
We claim that they are equivalent: Since the vertex groups $\varepsilon 
\mathcal{G}_{k\mid V}^{\prime p,p}$ and $\varepsilon \mathcal{G}_{k\mid
V}^{p,p}$ are conjugate inside $\mathcal{U}_{k+1}^{p,p}$ for all $p\in V,$
there exists some $a(p)\overset{def}{=}a_{k+1}^{p,p}\in \mathcal{U}%
_{k+1}^{p,p}$ with $a(p)\cdot \varepsilon \mathcal{G}_{k+1\mid
V}^{p,p}=\varepsilon ^{\prime }\mathcal{G}_{k+1\mid V}^{\prime p,p}$. Let $%
\mathcal{A}(p)$ denote the set of all such $a(p)$'s. The set $\mathcal{A}(p)$
is in 1-1 correspondence with the normalizer of $\mathcal{G}_{k+1}^{p,p}$
inside $\mathcal{U}_{k+1}^{p,p}$. We now have the bundle $\cup _{p\in U}%
\mathcal{A}(p)\rightarrow U$ which admits a crossection.

\QTP{Body Math}
We will conclude this section with three remarks.

\QTP{Body Math}
$1)$ We defined a Riemann geometry $\varepsilon \mathcal{G}_{1}$ in Example
3 using the pair $(g,\varepsilon )$ where $\varepsilon $ is the
LC-connection of $g$. Suppose $\varepsilon \mathcal{G}_{1}\sim \varepsilon
^{\prime }\mathcal{G}_{1}^{\prime }$ so that $\varepsilon ^{\prime }\mathcal{%
G}_{1}^{\prime }$ is defined by the geometric object $(g^{\prime
},\varepsilon ^{\prime }).$ Now $\varepsilon ^{\prime }$ need not be the LC
connection of $g^{\prime }!$ The reason is the Christoffel symbols $%
\varepsilon _{jk}^{i}$ can be expressed in terms of the derivatives of $%
g_{ij}$ whereas a gauge transformation is a section of jets and preserves
differentiation only pointwise at the level of jets but not locally. Now for 
$p\in M$, we can find a coordinate system around $p$ such that $%
g_{ij}=\delta _{ij}$ and $\varepsilon _{jk}^{i}=0$ at $p.$ The components $%
\varepsilon _{jk}^{\prime i}$ of $\varepsilon ^{\prime }$also satisfy this
condition at $p$ but we may not be able to express $\varepsilon
_{jk}^{\prime i}$ locally in terms of the derivatives of $g_{ij.\text{ }%
}^{\prime }$.

\QTP{Body Math}
Now there exists a somewhat stronger concept of equivalence of $phg$'s
which, in view of the proof of Proposition 34, exhibits a unique canonical
representative in each equivalence class which is "defined in terms of the
derivatives of some geometric object" and in particular gives the
Levi-Civita connection in Riemann geometry. This gives a generalization of
the main construction of [7] for parabolic geometries to all $phg$'s. The
idea is simple: for any $phg$ $\varepsilon \mathcal{G}_{k}$ we can define a
geometric object \textbf{g} of order $k+1$ on $M$ such that a $(k+1)$-arrow
of $\mathcal{U}_{k+1}$ belongs to $\varepsilon \mathcal{G}_{k}$ if and only
if it preserves \textbf{g}$.$ Therefore this condition gives the defining
equations of $\varepsilon \mathcal{G}_{k}.$ It is very easy to construct 
\textbf{g}$:$ For $x\in M$ consider the left coset space $\mathcal{U}%
_{k+1}^{x,x}/\varepsilon \mathcal{G}_{k}^{x,x}$ and define the bundle of
geometric objects $\mathcal{O}\overset{def}{=}\cup _{x\in M}\mathcal{U}%
_{k+1}^{x,x}/\varepsilon \mathcal{G}_{k}^{x,x}\rightarrow M$. Now $%
\varepsilon \mathcal{G}_{k}$ defines a global crossection of $\mathcal{%
O\rightarrow }M$ which is \textbf{g}$.$ In Example 3 \textbf{g}\emph{\ }$=$ $%
(g,\varepsilon )$ and in Example 4 \textbf{g}\emph{\ }$=$ $S=$ the
expression for the Schwarzian derivative! Now $a\in \Gamma \mathcal{A}_{k+1}$
acts on the sections of $\mathcal{O\rightarrow }M$ on the left. If $a\ast $%
\textbf{g }$=$ \textbf{g'} then \textbf{g'}\emph{\ }defines another $phg$ $%
a\ast \varepsilon \mathcal{G}_{k}$ which preserves \textbf{g' }and is
equivalent to $\varepsilon \mathcal{G}_{k}$ as defined by (59) because the
left coset $\beta \varepsilon \mathcal{G}_{k}^{x,x}$ defines the conjugate $%
\beta \varepsilon \mathcal{G}_{k}^{x,x}\beta ^{-1}$ but not conversely since
the normalizer of $\varepsilon \mathcal{G}_{k}^{x,x}$ inside $\mathcal{U}%
_{k+1}^{x,x}$ may strictly contain $\varepsilon \mathcal{G}_{k}^{x,x}$ in
general. In short, the philosophy of (59) is to preserve the symmetry group
of the object whereas the philosophy of the second is to preserve the object
itself.

\QTP{Body Math}
$2)$ Let $\Gamma \mathcal{A}_{k+1,i}$ denote the group of sections of $%
\mathcal{A}_{k+1,i}\rightarrow M$ defined as the kernel of the projection $%
\pi :\mathcal{A}_{k+1}\rightarrow \mathcal{A}_{i}$. If $a\in \Gamma \mathcal{%
A}_{k+1,i}$, then $\varepsilon \mathcal{G}_{k}$ and $a\cdot \varepsilon 
\mathcal{G}_{k}$ define the same principal bundle $\mathcal{G}_{i}^{e,\cdot
}\rightarrow M$ since $\pi (a)=Id.$ Thus we get the isomorphism of
algebroids $a^{\ast }:(\mathcal{G}_{i}^{e,\bullet },\varepsilon \mathfrak{G}%
_{k})\longrightarrow (\mathcal{G}_{i}^{e,\bullet },a\cdot \varepsilon 
\mathfrak{G}_{k})$ of order $k-i.$ In particular, let us choose $i=k.$
Recalling that $\varepsilon :\mathcal{G}_{k}\rightarrow \varepsilon \mathcal{%
G}_{k}$ is a connection on the principal bundle $\mathcal{G}_{k}^{e,\bullet
}\rightarrow M$, acting with $a\in \Gamma \mathcal{A}_{k+1,k}$ on $%
\varepsilon \mathcal{G}_{k}$ amounts to fixing $\mathcal{G}_{i}^{e,\bullet }$
but changing the connection. So the actions of $a\in \Gamma \mathcal{A}%
_{k+1,k}$ give points on the moduli space $\mathcal{M(G}_{k}^{e,\bullet })$
of connections on the principal bundle $\mathcal{G}_{k}^{e,\bullet
}\rightarrow M$ as defined in gauge theory whereas they keep us inside $%
\left[ \varepsilon \mathcal{G}_{k}\right] .$ It follows that $\left[
\varepsilon \mathcal{G}_{k}\right] $ plays the role of $\mathcal{M(G}%
_{k}^{e,\bullet })$ and (66) reflects the philosophy of attaching certain
invariants to $\mathcal{M(G}_{k}^{e,\bullet }).$

\QTP{Body Math}
$3)$ All the constructions in this note can be done in the holomorphic
category by considering $k$-jets of holomorphic objects and then working in
the smooth category as above. However many subtleties arise. For instance,
the standard definition of $\mathcal{P}^{\ast }(M,\varepsilon \mathfrak{G}%
_{k})$ in terms of the Chern classes of the complexification of $\varepsilon 
\mathfrak{G}_{k}\rightarrow M$ becomes rather artificial because the
underlying idea is the complexification of a homogeneous space which is a
nontrivial problem even for Lie groups (see [18], pg. 429-430 and [11]).

\section{Appendix A: Cartan connections}

We resume the setting of Section 7. Suppose $\mathcal{R}_{0}=0%
\Leftrightarrow \mathfrak{R}_{0}=0$ so that the $phg$ $\left( \mathcal{G}%
_{k}^{e,\bullet },\varepsilon \mathcal{G}_{k}\right) $ of order zero
integrates to the pseudogroup $\widetilde{\left( \mathcal{G}_{k}^{e,\bullet
},\varepsilon \mathcal{G}_{k}\right) }.$ Therefore $\widehat{\mathcal{P}%
^{\ast }}(\mathcal{G}_{k}^{e,\bullet },\varepsilon \mathfrak{G}_{k})=0.$ For
simplicity, assume that $\widetilde{\left( \mathcal{G}_{k}^{e,\bullet
},\varepsilon \mathcal{G}_{k}\right) }$ globalizes to a Lie group $G$ so
that $G$ acts simply transitively on $\mathcal{G}_{k}^{e,\bullet }$ and $%
\mathcal{G}_{k}^{e,\bullet }\simeq G.$ Recall that the action of $G$ may not
descend to $\mathcal{G}_{k-1}^{e,\bullet }$ since we may not have $\mathfrak{%
R}_{1}=0.$ Therefore, recalling (11), we see that the assumption $\mathcal{R}%
_{0}=0$ fixes also the group other than the vertex class in the definition
of the $phg$. Let $\mathfrak{g}$ denote the Lie algebra of the infinitesimal
generators of $G$ which can be identified with a subalgebra of $\mathfrak{X}(%
\mathcal{G}_{k}^{e,\bullet })$ = the Lie algebra of vector fields on $%
\mathcal{G}_{k}^{e,\bullet }.$ Recall that this identification is done by
evaluating the infinitesimal generators at some point and therefore is not
canonical.

Now the restriction of $D_{1}^{\prime }$ in (37) to the principal bundle $%
\mathcal{G}_{k}^{e,\bullet }\rightarrow M$ defines an $\mathfrak{h}$-valued $%
1$-form on $\mathcal{G}_{k}^{e,\bullet }$ where $\mathfrak{h}=$ the Lie
algebra of $\mathcal{G}_{k}^{e,e}.$ This is the connection in Proposition
18. If also $\mathfrak{R}_{0}=0$ as we now assume, we can define a $%
\mathfrak{g}$-valued $1$-form $\omega $ on $\mathcal{G}_{k}^{e,\bullet }$ as
follows. Let $\overline{p}=j_{k}(f)^{e,p}\in \mathcal{G}_{k}^{e,\bullet }$
and $\xi _{\overline{p}}\in T(\mathcal{G}_{k}^{e,\bullet })^{\overline{p}}=%
\mathfrak{G}_{k}^{p}\simeq \mathfrak{g}$. We rewrite (26) as

\begin{equation}
\varepsilon j_{k}(f^{-1})_{\ast }^{p,e}:T_{\overline{p}}(\mathcal{G}%
_{k}^{e,\bullet })^{\overline{p}}\rightarrow T(\mathcal{G}_{k}^{e,\bullet
})^{\overline{e}}=\mathfrak{g}
\end{equation}

We define $\omega \overset{def}{=}\varepsilon j_{k}(f^{-1})_{\ast }^{p,e}$
and easily show

\begin{proposition}
$\omega $ is a Cartan connection on $\mathcal{G}_{k}^{e,\bullet }\rightarrow
M.$
\end{proposition}

Let $\overline{\mathfrak{R}}$ denote the curvature of $\omega $ which is a $%
\mathfrak{g}$-valued $2$-form on $\mathcal{G}_{k}^{e,\bullet }.$ Since the
conditions $\overline{\mathfrak{R}}=0$ and $\mathfrak{R}_{k}=0$ are both
equivalent to local homogeneity of $M,$ obviously we have

\begin{proposition}
$\overline{\mathfrak{R}}=0\Leftrightarrow \mathfrak{R}_{k}=0$
\end{proposition}

Now $\mathfrak{R}_{k}(\overline{p})(\eta _{p},\sigma _{p}):\mathfrak{G}%
_{k}^{p}\rightarrow \mathfrak{G}_{k}^{p}$ becomes

\begin{equation}
\mathfrak{R}_{k}(\overline{p})(\eta _{p},\sigma _{p}):\mathfrak{g}%
\rightarrow \mathfrak{g}
\end{equation}

\begin{proposition}
(68) is a derivation.
\end{proposition}

Proposition 37 together with the interpretation of $Der(\mathfrak{g})$ in
[15] now gives a very interesting interpretation of $\mathfrak{R}_{k}.$

We recall the representation 
\begin{equation}
ad:\mathfrak{g\rightarrow }Der(\mathfrak{g})
\end{equation}
\ \ and assume

\textbf{A2:} (69) is an isomorphism.

For instance \textbf{A2} holds if $\mathfrak{g}$ is semisimple. However,
this assumption forces $k\leq 2$ and all our efforts with higher order jets
fall flat! Assuming \textbf{A2}, we identify $\mathfrak{R}_{k}(\overline{p}%
)(\eta _{p},\sigma _{p})$ uniquely with an element of $\mathfrak{g}$ so that 
$\mathfrak{R}_{k}$ becomes a $\mathfrak{g}$-valued $2$-form on $\mathcal{G}%
_{k}^{e,\bullet }\rightarrow M.$

\begin{proposition}
If \textbf{A2} holds, then $\overline{\mathfrak{R}}=\mathfrak{R}_{k}.$
\end{proposition}

Thus we conclude that $\mathfrak{R}_{0}$ and therefore $\widehat{\mathcal{P}%
^{\ast }}(\mathcal{G}_{k}^{e,\bullet },\varepsilon \mathfrak{G}_{k})$ is an
obstruction to the existence of a Cartan connection on the principal bundles 
$a\cdot \mathcal{G}_{k}^{e,\bullet }\rightarrow M$, $a\in \Gamma \mathcal{A}%
_{k+1}.$

Assuming $\mathfrak{R}_{0}=0,$ now $\mathfrak{R}_{1}$ and therefore $%
\widehat{\mathcal{P}^{\ast }}(\mathcal{G}_{k-1}^{e,\bullet },\varepsilon 
\mathfrak{G}_{k})$ is an obstruction to the existence of a (generalized)
Cartan connection on the principal bundles $\pi _{k,k-1}(a)\cdot \mathcal{G}%
_{k-1}^{e,\bullet }\rightarrow M$, $a\in \Gamma \mathcal{A}_{k+1}$ in a way
which is straightforward at this stage. With the assumptions $\mathcal{R}%
_{i}=0$, $0\leq i\leq k-1,$ we finally encounter the topological
obstructions $\widehat{\mathcal{P}^{\ast }}(\mathcal{G}_{0}^{e,\bullet
},\varepsilon \mathfrak{G}_{k})$ to the homogeneity of $M$ by the action of $%
G.$

\section{Appendix B: Chern-Simons forms}

In this section $\varepsilon \mathcal{G}_{k}$ is a $phg$ on $M$ with $k=0.$
Equivalently, we have a splitting $\varepsilon :M\times M\rightarrow 
\mathcal{U}_{1}$ which in turn is equivalent to the parallelizability of $M.$
We refer to [1] for a detailed study of this case. According to Proposition
24, the total space of the principal bundle $\mathcal{G}_{k}^{e,\bullet
}\rightarrow M$ defined by $\varepsilon \mathcal{G}_{k}$ is parallelizable
which gives a rich source of examples.

Now since $M$ is parallelized by $\varepsilon $, $\widehat{\mathcal{P}^{\ast
}}(\mathcal{G}_{0}^{e,\bullet },\varepsilon \mathfrak{G}_{0})=$ $\mathcal{P}%
^{\ast }(M,T)=0$ and therefore the forms $\mathfrak{R}_{0}^{2i}\in \wedge
^{4i}T^{\ast }$ defined by (52) are exact. Our purpose here is to show that
the "Chern-Simons" forms (but with a surprisingly different interpretation)
furnish some canonical primitives of these forms. Henceforth we denote the
curvature $\mathfrak{R}_{0}$ by $\mathfrak{R.}$

First, we recall from [1] the definition of the curvature $\widetilde{%
\mathfrak{R}}$ defined by

\begin{equation}
\widetilde{\mathfrak{R}}_{rj,k}^{i}=\left[ \frac{\partial \Gamma _{jk}^{i}}{%
\partial x^{r}}+\Gamma _{rk}^{a}\Gamma _{ja}^{i}\right] _{[rj]}
\end{equation}%
where

\begin{equation}
\Gamma _{jk}^{i}(x)\overset{def}{=}\left[ \frac{\partial \varepsilon
_{k}^{i}(x,y)}{\partial y^{j}}\right] _{y=x}
\end{equation}

We always have $\widetilde{\mathfrak{R}}=0$ on a parallelizable manifold $%
(M,\varepsilon ).$ The reason is that $\widetilde{\mathfrak{R}}=0$ gives the
integrability conditions of

\begin{equation}
\widetilde{\nabla }_{r}\xi ^{i}\overset{def}{=}\frac{\partial \xi ^{i}}{%
\partial x^{r}}-\Gamma _{ra}^{i}\xi ^{a}=0
\end{equation}%
and a vector field $\xi =(\xi ^{i})$ solves (72) if ond only if it is $%
\varepsilon $-invariant. Since we start with the global parallelism $%
\varepsilon ,$ we can always construct $\varepsilon $-invariant vector
fields with arbitrary initial conditions and therefore $\widetilde{\mathfrak{%
R}}=0.$ However

\begin{equation}
\mathfrak{R}_{rj,k}^{i}=\left[ \frac{\partial \Gamma _{kj}^{i}}{\partial
x^{r}}+\Gamma _{kr}^{a}\Gamma _{aj}^{i}\right] _{[rj]}
\end{equation}%
and $\mathfrak{R}=0$ if and only if $M$ is locally homogeneous in which case 
$(M,\varepsilon )$ is called a local Lie group in [1].

Using $\widetilde{\mathfrak{R}}=0,$ we will now construct a locally exact
complex. Consider the vector bundle $T^{\ast }\otimes T\rightarrow M$
isomorphic to $Hom(T,T)\rightarrow M$ and the vector bundle $\wedge
^{k}T^{\ast }\otimes Hom(T,T)\rightarrow M$, the bundle of $k$-forms on $M$
with values in $Hom(T,T).$ A local section of $\wedge ^{k}T^{\ast }\otimes
Hom(T,T)\rightarrow M$ is of the form $\omega _{h_{k}....h_{1},j}^{i}$ where 
$\omega $ is alternating in the indices $h_{k},...,h_{1}.$ We define the
local operator $\widetilde{d}_{r}$ by the formula

\begin{equation}
\widetilde{d}_{r}\omega _{h_{k}....h_{1},j}^{i}\overset{def}{=}\frac{%
\partial \omega _{h_{k}....h_{1},j}^{i}}{\partial x^{r}}-\Gamma
_{ra}^{i}\omega _{h_{k}....h_{1},j}^{a}+\Gamma _{rj}^{a}\omega
_{h_{k}....h_{1},a}^{i}
\end{equation}

The operator $\widetilde{d}_{r}$ has a coordinatefree meaning only for $k=0$
in which case $\widetilde{d}_{r}=\widetilde{\nabla }_{r}$ which is defined
as an extension of (72) on arbitrary tensor fields. Now we define the first
order linear differential operator

\begin{equation}
\widetilde{d}:\wedge ^{k}T^{\ast }\otimes Hom(T,T)\longrightarrow \wedge
^{k+1}T^{\ast }\otimes Hom(T,T)
\end{equation}%
by the formula

\begin{eqnarray}
&&\left( \widetilde{d}\omega \right) _{rh_{k}...h_{1},j}^{i}\overset{def}{=}%
\left[ \widetilde{d}_{r}\omega _{h_{k}....h_{1},j}^{i}\right]
_{[rh_{k}...h_{1}]} \\
&=&\widetilde{d}_{r}\omega _{h_{k}....h_{1},j}^{i}-\widetilde{d}%
_{h_{k}}\omega _{rh_{k-1}....h_{1},j}^{i}-...-\widetilde{d}_{h_{1}}\omega
_{h_{k}....h_{2}r,j}^{i}  \notag
\end{eqnarray}

Since $\widetilde{\mathfrak{R}}=0$, we have $\widetilde{d}\circ \widetilde{d}%
=0$ and we obtain the complex

\begin{equation}
Hom(T,T)\overset{\widetilde{d}}{\longrightarrow }\wedge ^{1}T^{\ast }\otimes
Hom(T,T)\overset{\widetilde{d}}{\longrightarrow }....\overset{\widetilde{d}}{%
\longrightarrow }\wedge ^{n}T^{\ast }\otimes Hom(T,T)
\end{equation}%
which is locally exact. It is easy to give a coordinate free description of
(77) which is a well known construction. However, observe that (77) is not a
complex with the accordingly defined operators $d_{r}$ since we do not
assume $\mathfrak{R}=0$ (see (93) below)$.$ The kernel $\widetilde{Hom(T,T)}$
of the first operator in (77) is $\varepsilon $-invariant sections of $%
Hom(T,T)\rightarrow M$ and (77) is a fine resolution of the sheaf $%
\widetilde{Hom(T,T)}.$

Now let $\omega \in \wedge ^{k}T^{\ast }\otimes Hom(T,T)$ and $\psi \in
\wedge ^{m}T^{\ast }\otimes Hom(T,T).$ We define $\omega \wedge \psi $ by
the formula

\begin{equation}
(\omega \wedge \psi )_{h_{k}...h_{1}s_{m}...s_{1},j}^{i}\overset{def}{=}%
\left[ \omega _{h_{k}...h_{1},a}^{i}\psi _{s_{m}...s_{1},j}^{a}\right]
_{[h_{k}...h_{1}s_{m}...s_{1}]}
\end{equation}%
or in coordinatefree form $\left( \omega \wedge \psi \right)
(X_{1},...X_{k},X_{k+1},...,Y_{k+m})\overset{def}{=}$

\begin{equation}
\dsum\limits_{\sigma \in S_{k+m}}\frac{1}{k!m!}sgn(\sigma )\omega (X_{\sigma
(1)},...,X_{\sigma (k)})\circ \psi (X_{\sigma (k+1)},...,Y_{\sigma (k+m})
\end{equation}%
where $\circ $ denotes composition in $Hom(T,T).$ We have

\begin{equation}
\widetilde{d}(\omega \wedge \psi )=\left( \widetilde{d}\omega \right) \wedge
\psi +(-1)^{\deg (\omega )}\omega \wedge \left( \widetilde{d}\psi \right)
\end{equation}

We recall the definition

\begin{equation}
T_{k,j}^{i}\overset{def}{=}\Gamma _{kj}^{i}-\Gamma _{jk}^{i}
\end{equation}

Now $T=(T_{k,j}^{i})\in $ $\wedge ^{1}T^{\ast }\otimes Hom(T,T)$ and $%
\mathfrak{R=}(\mathfrak{R}_{km,j}^{i})\in \wedge ^{2}T^{\ast }\otimes
Hom(T,T).$

We have the following fundamental

\begin{proposition}
(Structure equation)%
\begin{equation}
\widetilde{d}T+T\wedge T=\mathfrak{R}
\end{equation}
\end{proposition}

It is worthwhile to observe the remarkable analogy between (82) and the well
known structure equation

\begin{equation}
dA+A\wedge A=R
\end{equation}%
on a principal bundle $P\rightarrow M$ where $d$ is the exterior derivative, 
$A$ the Lie algebra valued connection $1$-form and $R$ its curvature $2$%
-form. Observe that (82) is defined on $M$ whereas (83) is defined on $P.$
To make this analogy more precise, let $P\rightarrow M$ be $\mathcal{U}%
_{1}^{e,\bullet }\rightarrow M$ with structure group $\mathcal{U}%
_{1}^{e,\bullet }\simeq GL(n,\mathbb{R})$ with Lie algebra $\mathfrak{h=gl}%
(n,\mathbb{R}).$ Now $\Gamma _{jk}^{i}$ defined in terms of the absolute
parallelism $\varepsilon $ by the formula (71) transform as the components
of $\mathfrak{gl}(n,\mathbb{R})$-valued $1$-form on $\mathcal{U}%
_{1}^{e,\bullet }\rightarrow M$ with curvature $R.$ In this particular case $%
Hom(T,T)\simeq \mathfrak{gl}(n,\mathbb{R})$ (see the last paragraph of
Section 4) and therefore $R$ and $\mathfrak{R}$ live in the same space...but
they are again different because $R=\widetilde{\mathfrak{R}}=0$ whereas $%
\mathfrak{R}$ need not vanish.

Proposition 39 shows that $\mathfrak{R}$ is determined by $T.$ In fact, we
have the following fundamental

\begin{proposition}
\begin{equation}
\widetilde{\nabla }_{r}T_{jk}^{i}=\mathfrak{R}_{jk,r}^{i}
\end{equation}
\end{proposition}

Therefore $\mathfrak{R}=0$ $\Leftrightarrow T$ is $\varepsilon $-invariant.
Observe that $\mathfrak{R}$ and $T$ have the same alternating indices $j,k$.

Now let $\omega \in \wedge ^{k}T^{\ast }\otimes Hom(T,T).$ We define $%
Tr\left( \omega \right) \in \wedge ^{k}T^{\ast }$ by the formula

\begin{equation}
Tr(\omega _{h_{k}...h_{1},j}^{i})\overset{def}{=}\omega
_{h_{k}...h_{1},a}^{a}
\end{equation}

So we obtain the following commutative diagram

\begin{equation}
\begin{array}{ccccccc}
Hom(T,T) & \overset{\widetilde{d}}{\longrightarrow } & \wedge ^{1}T^{\ast
}\otimes Hom(T,T) & \overset{\widetilde{d}}{\longrightarrow } & .... & 
\overset{\widetilde{d}}{\longrightarrow } & \wedge ^{n}T^{\ast }\otimes
Hom(T,T) \\ 
\downarrow Tr &  & \downarrow Tr &  & \downarrow Tr &  & \downarrow Tr \\ 
C^{\infty }(M) & \overset{d}{\longrightarrow } & \wedge ^{1}T^{\ast } & 
\overset{d}{\longrightarrow } & .... & \overset{d}{\longrightarrow } & 
\wedge ^{n}T^{\ast }%
\end{array}%
\end{equation}%
where the lower complex in (86) is the de Rham complex of $M.$

With the notation (52), we now have 
\begin{equation}
\mathfrak{R}^{i}=\mathfrak{R\wedge R\wedge ....\wedge R}\text{ \ \ }(i\text{%
-copies)}
\end{equation}%
and we define

\begin{equation}
T^{i}\overset{def}{=}T\wedge T\wedge ...\wedge T\text{ \ \ \ }(i\text{%
-copies)}
\end{equation}

Clearly $\mathfrak{R}^{i}\in \wedge ^{2i}T^{\ast }\otimes Hom(T,T)$ and $%
T^{i}\in \wedge ^{i}T^{\ast }\otimes Hom(T,T)$. Therefore $Tr(\mathfrak{R}%
^{i})\in \wedge ^{2i}T^{\ast }$ and $Tr(T^{i})\in \wedge ^{i}T^{\ast }.$ It
is easy to see that

\begin{equation}
Tr(T^{2i})=0
\end{equation}%
and we are left with $Tr(T^{2i+1}),$ $i\geq 0.$ Omitting $\wedge $ from our
notation, applying $\widetilde{d}$ to (82) and substituting back from (82)
we obtain

\begin{equation}
\widetilde{d}\mathfrak{R}=\mathfrak{R}T-T\mathfrak{R}
\end{equation}%
Using (82), (90) and (80) we compute

\begin{eqnarray}
\widetilde{d}(T^{3}) &=&\mathfrak{R}T^{2}-T\mathfrak{R}T+T^{2}\mathfrak{R-}%
T^{4} \\
\widetilde{d}(\mathfrak{R}T) &=&-T\mathfrak{R}T+\mathfrak{R}^{2}  \notag
\end{eqnarray}%
Taking the trace of the formulas in (91) and observing $Tr(T\mathfrak{R}%
T)=-Tr(\mathfrak{R}T^{2})=-Tr(T^{2}\mathfrak{R}),$ we deduce

\begin{equation}
dTr\left( \mathfrak{R}T-\frac{1}{3}T^{3}\right) =Tr(\mathfrak{R}^{2})
\end{equation}

Observe the "Chern-Simons" $3$-form in (92) with the surprising difference
that the Lie algebra valued $1$-form $A$ in (83) is replaced with the $%
Hom(T,T)$-valued $1$-form $T$ and it lives on the base $M$!! The higher
degree Chern-Simons forms are derived in the same way without any further
computation but as a logical consequence of the correspondence between (82)
and (83).

To complete the analogy to the formalism of connections on principal
bundles, we recall the Bianchi identity

\begin{equation}
DR=0
\end{equation}%
on $P\rightarrow M$ where $D$ is the exterior covariant differentiation. We
define

\begin{equation}
d_{r}\omega _{h_{k}....h_{1},j}^{i}\overset{def}{=}\frac{\partial \omega
_{h_{k}....h_{1},j}^{i}}{\partial x^{r}}-\Gamma _{ar}^{i}\omega
_{h_{k}....h_{1},j}^{a}+\Gamma _{jr}^{a}\omega _{h_{k}....h_{1},a}^{i}
\end{equation}%
and the operator

\begin{equation}
d:\wedge ^{k}T^{\ast }\otimes Hom(T,T)\longrightarrow \wedge ^{k+1}T^{\ast
}\otimes Hom(T,T)
\end{equation}%
by the formula (76) using $d_{r}$ instead of $\widetilde{d}_{r}.$ Observe
that

\begin{equation}
d_{r}\omega _{h_{k}....h_{1},j}^{i}=\widetilde{d}_{r}\omega
_{h_{k}....h_{1},j}^{i}-T_{ar}^{i}\omega
_{h_{k}....h_{1},j}^{a}+T_{rj}^{a}\omega _{h_{k}....h_{1},a}^{i}
\end{equation}%
Now $d\circ d\neq 0$ because $\mathfrak{R}$ is the obstruction to $d\circ
d=0.$

The analog of (83) is

\begin{proposition}
(Bianchi identity) $d\mathfrak{R}=0$
\end{proposition}

Finally, the first formula in (91) shows $dTr(T^{3})=0$ if $\mathfrak{R}=0$.
An easy induction gives

\begin{proposition}
If $\mathfrak{R}=0$, then $dTr(T^{2i+1})=0.$

\begin{definition}
$\left[ Tr(T^{2i+1})\right] \in H_{dR}^{2i+1}(M,\mathbb{R})$ are the
secondary characteristic classes of the local Lie group $(M,\varepsilon ).$
\end{definition}
\end{proposition}

The secondary characteristic classes coincide with the Chern-Simons classes
on a local Lie group.

We will conclude with a question. Recall that the forms $Tr(\mathfrak{R}%
_{k-i}^{j})\in $ $H_{dR}^{2j}(\varepsilon \mathcal{G}_{i}^{e,\bullet },$ $%
\mathbb{R})$ are exact for $1\leq i\leq k$ and for $i=k$ we found above some
explicit primitives as "Chern-Simons" forms.

$\mathbf{Q:}$ Find some explicit primitives for $1\leq i\leq k-1.$

\section{Appendix C. Uniformization number and representations}

Let $(\mathfrak{g,h)}$ be a Lie algebra pair with $\mathfrak{h}\subset 
\mathfrak{g.}$ We set $V=\mathfrak{g/h}$ and define $ad_{\mathfrak{h,g/h}}:%
\mathfrak{h}\rightarrow gl(V)$ by

\begin{equation}
ad_{\mathfrak{h,g/h}}(h)(g+\mathfrak{h)}\overset{def}{=}[h,g]+\mathfrak{h}
\end{equation}

Clearly $ad_{\mathfrak{h,g/h}}$ is well defined and is a representation of $%
\mathfrak{h}$.

\begin{definition}
$ad_{\mathfrak{h,g/h}}$ is the adjoint representation of $\mathfrak{h}$
relative to $\mathfrak{g/h.}$
\end{definition}

Suppose $(\mathfrak{g}^{\prime }\mathfrak{,h)}$ is another such pair. We
call $(\mathfrak{g,h)}$ and $(\mathfrak{g}^{\prime }\mathfrak{,h)}$
isomorphic and write $(\mathfrak{g,h)\simeq }$ $(\mathfrak{g}^{\prime }%
\mathfrak{,h)}$ if $\mathfrak{g\simeq g}^{\prime }$ and the isomorphism $%
\simeq $ restricts to identity on $\mathfrak{h.}$ We easily check that if $(%
\mathfrak{g,h)\simeq }$ $(\mathfrak{g}^{\prime }\mathfrak{,h)}$ then the
representations $ad_{\mathfrak{h,g/h}}$ and $ad_{\mathfrak{h,g}^{\prime }/%
\mathfrak{h}}$ are isomorphic but not conversely.

\begin{proposition}
Any representation $\rho :\mathfrak{h\rightarrow }gl(W)$ is an adjoint
representation relative to some $\mathfrak{g\supset h}$.
\end{proposition}

Indeed, given a representation $\rho :\mathfrak{h\rightarrow }gl(W)$, we set 
$\mathfrak{g}\overset{def}{=}\mathfrak{h\times }W$ and check that $\mathfrak{%
g}$ is a Lie algebra with the bracket defined by $\left[ (h,w),(h^{\prime
},w^{\prime })\right] \overset{def}{=}([h,h^{\prime }],\rho (h)$ $w^{\prime
}-\rho (h^{\prime })w).$ We identify $\mathfrak{h}$ with the subalgebra $(%
\mathfrak{h},0)\subset \mathfrak{g}$ and $W$ with $\mathfrak{g}/\mathfrak{h}$
and check that $\rho =ad_{\mathfrak{h,g/h}}$ with these identifications.

Therefore adjoint representations exhaust all representations! From the
above construction of the pair $(\mathfrak{h\times }W,\mathfrak{h)}$ we
deduce

\begin{proposition}
Let $\mathfrak{h}$ be a (finite dimensional) Lie algebra. The following are
equivalent.
\end{proposition}

$i)$ (Ado's Theorem) $\mathfrak{h}$ has a faithful representation

$ii)$ There exists a (finite dimensional) Lie algebra $\mathfrak{g\supset h}$
such that the pair $(\mathfrak{g},\mathfrak{h})$ is effective.

Indeed, if $\rho :\mathfrak{h\rightarrow }gl(W)$ is faithful, then the above
pair $(\mathfrak{h\times }W,\mathfrak{h)}$ is effective and conversely, for
an effective pair $(\mathfrak{g},\mathfrak{h})$ the kernel of $Ad_{\mathfrak{%
h,g}}:\mathfrak{h\rightarrow }gl(\mathfrak{g})$ is $Z(\mathfrak{g})\cap 
\mathfrak{h}=\{0\}.$ Observe that $ord(\mathfrak{h\times }W,\mathfrak{h)=}1$%
. Recall that a Lie algebra $\mathfrak{g\supset h}$ with $(\mathfrak{g},%
\mathfrak{h})$ effective defines a "flag" inside $Nil(\mathfrak{h})$
according to (13). The next proposition therefore gives a far reaching
generalization of the Ado's theorem.

\begin{proposition}
Let $\mathfrak{h}$ be a (finite dimensional) Lie algebra with a flag $%
\mathcal{F}$ inside $Nil(\mathfrak{h})$. Then there exists a Lie algebra $%
\mathfrak{g\supset h}$ such that $(\mathfrak{g},\mathfrak{h})$ is effective
and defines $\mathcal{F}.$ In particular, there exists an effective pair $(%
\mathfrak{g},\mathfrak{h})$ with $ord(\mathfrak{g},\mathfrak{h})=\dim Nil(%
\mathfrak{h}).$
\end{proposition}

Finally, let $Iso_{k}(\mathfrak{h})$ denote the set of the above isomorphism
classes with $\dim \mathfrak{g-}\dim \mathfrak{h}$ $=k\geq 1.$ This set is
obtained from the set of \textit{all }representations of $\mathfrak{h}$ of
rank $k$ but the concept of isomorphism is more stringent than the usual
concept of isomorphism of two representations. For the Lie algebras $%
\mathfrak{g}_{i},$ $i=1,2,3$ in Example 3, for instance, the representations 
$\rho _{i}:\mathfrak{o}(n)\rightarrow \mathfrak{g}_{i}/\mathfrak{o}(n)$ are
isomorphic whereas the pairs $(\mathfrak{g}_{i},\mathfrak{o}(n))$ are
mutually nonisomorphic. Clearly the cardinality $\sharp Iso_{k}(\mathfrak{h}%
) $ is equal to the uniformization number $\sharp (\mathfrak{h,}k)$ defined
in Section 2 if we assume effectiveness.

\section{Appendix D. The adjoint representation}

Let $(M,\varepsilon \mathcal{G}_{k})$ be a locally solvable $phg$ and recall
the presheaf $\mathfrak{g}(U)$ whose sections are the local solutions of $%
\varepsilon \mathfrak{G}_{k}\rightarrow M$ on $U.$ We choose some $e\in U$
and define

\begin{eqnarray}
j_{k}(\cdot )^{e,e} &:&\mathfrak{g}(U)\rightarrow \left( \mathfrak{G}%
_{k}\right) ^{e} \\
&:&\xi \rightarrow j_{k}(\xi )^{e,e}  \notag
\end{eqnarray}

Now $\left( \mathfrak{G}_{k}\right) ^{e}$ is a Lie algebra endowed with the
algebraic bracket (19) and (98) is an isomorphism of Lie algebras for
sufficiently small and simply connected $U.$ Let $H^{\ast }(M,\mathfrak{g})$
denote the cohomology groups of $M$ with coefficients in the sheaf $%
\mathfrak{g.}$ Since $\mathfrak{g}$ is the kernel of the first operator in
(45) and partition of unity applies to sections of the spaces in (45), (45)
is a fine resolution of the sheaf $\mathfrak{g}$ and therefore computes $%
H^{\ast }(M,\mathfrak{g}).$ For simplicity of notation, we denote the Lie
algebra $\left( \mathfrak{G}_{k}\right) ^{e}$ by $\mathfrak{g}_{e}$ and let $%
H^{\ast }(\mathfrak{g}_{e}\mathfrak{,g}_{e})$ denote the deformation
cohomology of $\mathfrak{g}_{e}.$

\begin{proposition}
If $M$ is compact and simply connected, then 
\begin{equation}
H^{\ast }(M,\mathfrak{g})\simeq H^{\ast }(\mathfrak{g}_{e}\mathfrak{,g}_{e})
\end{equation}
\end{proposition}

For simplicity we now assume that the pseudogroup $\widetilde{(M,\varepsilon 
\mathcal{G}_{k})}$ is globalizable to $G$ so that any initial condition of $%
\left( \mathfrak{G}_{k}\right) ^{e}$ in (98) comes from a global section of $%
\mathfrak{g}(M)=$ the Lie algebra of infinitesimal generators of the action
of $G$ on $M.$ Now $G$ acts on $\Gamma \mathfrak{G}_{k}=$ the space of the
global sections of the algebroid $\mathfrak{G}_{k}\rightarrow M$ by

\begin{equation}
\left( g\cdot s\right) (x)\overset{def}{=}\varepsilon j_{k}(g)_{\ast
}^{x,g(x)}s(g^{-1}(x))
\end{equation}

Thus $\Gamma \mathfrak{G}_{k}$ is an infinite dimensional representation
space for $G$ and $\mathfrak{g}(M)\subset \Gamma \mathfrak{G}_{k}$ is a
stable and finite dimensional subspace. The adjoint representation of $G$ on 
$\mathfrak{g}(M)$ localizes as follows: We identify $\xi \in \mathfrak{g}(M)$
with $j_{k}(\xi )^{e,e}$ by (98). Now $j_{k}(\xi )^{g(e),g(e)}$ determines
some $\zeta \in \mathfrak{g}(M)$ and $Ad(g)^{e}$ is the map $j_{k}(\xi
)^{e,e}\rightarrow j_{k}(\zeta )^{e,e}.$

\bigskip

\textbf{References}

\bigskip

[1] E.Abado\u{g}lu, E.Orta\c{c}gil: Intrinsic characteristic classes of a
local Lie group, Portugal. Math. (N.S.), Vol.67, Fasc.4 (2010), 453-483

[2] E.Abado\u{g}lu, E.Orta\c{c}gil, F. \"{O}zt\"{u}rk: Klein geometries,
parabolic geometries differential equations of finite type, J. Lie Theory 18
(2008), 67-82

[3] A.D.Blaom: Geometric structures as deformed infinitesimal symmetries,
Trans. Amer. Soc. 358 (2006), 3651-3671

[4] A.D.Blaom: Lie algebroids and Cartan's method of equivalence, Trans.
Amer. Soc. 364 (2012), 3071-3135

[5] A.D.Blaom: The infinitesimalization and reconstruction of locally
homogeneous manifolds, to appear in SIGMA

[6] A.D.Blaom: A Lie theory of pseudogroups, preprint

[7] R.Bott: Lectures on characteristic classes and foliations, Springer
Lecture Notes, 279, 1972

[8] A.Cap, A.R.Gover: Tractor calculi for parabolic geometries, Trans. Amer.
Math. Soc., 2002, 1511-1548

[9] G.A.Fredricks, P.B.Gilkey, P.E.Parker: A higher order invariant of
differential manifolds, Trans. Amer. Soc., 315, 1989, 373-388 \ 

[10] V.M.Goldman: Locally homogeneous geometric manifolds, In Proceedings of
the International Congress of Mathematics, Vol. II, Hindustan Book Agency,
New Delhi, 2010, 717-744

[11] G.Hochschild: Complexification of real analytic groups, Trans. Amer.
Math. Soc., 125, 1966, 406-413

[12] A.Kumpera, D.C.Spencer: Lie equations, Vol.1, General Theory, Ann.
Math. Stud. 73, 1972

[13] D.Montgomery: Simply connected homogeneous spaces, Proceedings of the
A.M.S, Vol.1, (1950), 467-469

[14] P.Olver, J.Pohjanpelto: Maurer Cartan forms and the structure of Lie
pseudogroups, Selecta Math. 11, 2005, 99-126

[15] P.Olver, E.Orta\c{c}gil, M.Ta\c{s}k\i n: The symmetry group of a local
Lie group, in progress

[16] J.P.Pommaret: Systems of partial differential equations and Lie pseudo
groups, Math. Appl.. 14, Gordon \& Breach Science Publishers, 1978

[17] J.F.Pommaret: Partial differential equations and group theory, Math.
Appl. 293, Kluwer Academic Publishers, 1994

[18] L.Pontryagin: Topological groups, Selected Works, Vol.2, Gordon \&
Breach Science Publishers, 1986

[19] G.Weingart: Holonomic and semi-holonomic geometries, In Global analysis
and harmonic analysis, Marseille-Luminy, Semin.Congr. 4, Soc. Math. France,
Paris, 2000, 307-328

\bigskip

\bigskip

Erc\"{u}ment Orta\c{c}gil, Bodrum, T\"{u}rkiye

ortacgile@gmail.com

ortacgil@boun.edu.tr

\end{document}